\documentclass[12pt]{article}
\usepackage{amsfonts}
\usepackage{amsmath}
\usepackage[english]{babel}
\usepackage{amsthm}
\usepackage{tikz-cd}
\usepackage{mathtools}
\usepackage[all,cmtip]{xy}
\usepackage{mathrsfs}
\usepackage[all,cmtip]{xy}
\usepackage{graphicx}
\usepackage{titlesec} 
\usepackage[T1]{fontenc}
\usepackage{newtxtext}
\usepackage{newtxmath}
\usepackage{indentfirst}
\usepackage{biblatex} 
\usepackage{enumitem}
\usepackage[colorlinks=true, linkcolor=blue, citecolor=blue, urlcolor=blue]{hyperref}
\usepackage[a4paper,left=1 in,right=1 in , bottom = 1.5 in]{geometry}
\numberwithin{equation}{section}
\theoremstyle{definition}
\newtheorem{definition}{Definition}[section]
\newtheorem{exam}{Example}
\newtheorem{Thm}{Theorem}[section]
\newtheorem{Lem}[Thm]{Lemma}
\newtheorem{Prpn}[Thm]{Proposition}
\newtheorem{rmk}{Remark}
\newtheorem{Cor}[Thm]{Corollary}
\titleformat{\section}{\normalfont\Large\bfseries}{\S\thesection}{1em}{}[] 

\title{Arithmetic jet stacks}
\author{Akhil Surendran}
\date{}

\begin{document}
	\maketitle
	\begin{abstract}
		We construct a stack, which generalizes the notion of the arithmetic jet space functor constructed by Borger using the generalized Witt vector constructions. A priori we construct the jet space functor for prestacks, consequently showing that this functor preserves stacks. We then proceed to show the representability of the jet space functor applied to a representable morphism of stacks thereby showing that the jet space functor preserves algebraic and Deligne-Mumford Stacks. We also discuss the properties of morphisms preserved under the jet functor.
	\end{abstract}
	
	\section*{Introduction}
	
	 Differential algebraic methods were developed by Buium to study Diophantine problems through techniques inspired by differential geometry \cite{article5}. In the arithmetic setting, the absence of non-trivial derivations on number fields prevents a direct application of classical differential algebra. To overcome this obstacle, Buium introduced the notion of a $p$-derivation, an arithmetic analogue of an ordinary derivation \cite{article4}. Correspondingly, developed the theory of arithmetic jet spaces, which serve as arithmetic analogues of classical jet spaces (see Section \ref{sec:Preliminaries} for the more general notion of $\pi$-derivations). The theory of arithmetic jet spaces has found several remarkable applications in arithmetic geometry. Using first-order $p$-jet spaces, Buium obtained a proof of the Manin--Mumford conjecture together with an effective bound \cite{article4}. The study of $p$-jet spaces associated with the Hodge bundle on modular curves led to the introduction of $\delta$-modular forms \cite{article1}. Furthermore, Buium's theory of $\delta$-characters \cite{article6} was used by Borger and Saha in \cite{article8} to construct a filtered $F$-isocrystal associated with an abelian scheme over a $p$-adically complete discrete valuation ring with perfect residue field. This construction was subsequently extended by Pandit and Saha to the setting of $\pi$-formal group schemes \cite{article7}. Gurney and Buium extended the Serre–Tate theory of canonical lifts of ordinary abelian varieties to arbitrary unpolarised families of ordinary abelian varieties parameterised by a p-adic formal scheme \cite{article14}. More recently, Buium and Miller established a connection between $\delta$-geometry and perfectoid geometry by associating perfectoid spaces to smooth schemes and morphisms of perfectoid spaces to $\delta$-morphisms of smooth schemes \cite{article10}.
	 
   	 Borger constructed the generalized Witt vector functor $W_{R,E,n}$ of length $n$ in \cite{article2} from the category of $R$-algebras to itself; see \cite[Section~1]{article2} or Section \ref{sec:Preliminaries} for the construction and assumptions on $R$. The generalized Witt vector functor retains several key properties previously observed in the p-typical case such as the preservation of étale morphisms and the existence of a left adjoint given by $\Lambda _{R,E}$, which reduces to the usual $p$-jet spaces in the $p$-typical case. In \cite{article3}, Borger generalized these functors to the category of sheaves of sets. The construction is as follows; let \( \text{Aff}_S \) denote the category of affine schemes over $S$, where $S$ is a separated scheme and further assumptions on it are given in Section \ref{sec: COnstruction of jets} which is fixed throughout this paper. Let $X$ be a sheaf of sets over \( \text{Aff}_S \) equipped with étale topology. Borger's generalized Witt vector functor $W_{R,E,n}$ induces a functor on sheaves over the affine étale site. More precisely, if $X$ is a sheaf on $\text{Aff}_S$, we define
	\begin{center}
		$W_{S,E,n*}X := X\circ W_{R,E,n}.$
	\end{center}
	As shown in \cite{article3}, this construction preserves sheaves and hence determines an endofunctor on $ Sp_S $, where $Sp_S$ denotes the category of sheaves on $\text{Aff}_S$. Using the equivalence between the categories of sheaves on $Sch/S$ and $\text{Aff}_S$ , we freely regard $W_{S,E,n*}$ as an endofunctor on $Sp_S$ over $Sch/S$. The functor $W_{S,E,n*}$ has a left adjoint $W_{S,E,n}^{*}$ which extends the generalized Witt vector functor in the level of rings to sheaves. For the rest of this section, we denote $W_{S,E,n*}$ as $W_{n*}$.
	
	Let $\mathscr{X}$ be a prestack (according to \cite{textbook2}) over $Sch/S$, we construct a category $J^{n}\mathscr{X}$ along with a functor $p:J^{n}\mathscr{X}\rightarrow Sch/S$ such that the data $(J^{n}\mathscr{X},p)$ satisfies the following
	
	\begin{Thm}
		If $\mathscr{X}$ is a stack then the prestack $J^{n}\mathscr{X}$ is a stack.
	\end{Thm}
	Let $G$ be a group scheme over $S$ acting on an $S$-scheme $U$. In \cite{article3}, Borger showed that the functor $W_{n*}$ preserves schemes. Moreover, the action of $G$ on $U$ induces, via the functoriality of $W_{n*}$, an action of $W_{n*}G$ on $W_{n*}U$. In particular, if $G$ is a group scheme, then $W_{n*}G$ naturally inherits the structure of a group scheme. Next result shows that the stackification of arithmetic jets of the prestack obtained by $G$-action of $U$ is the quotient stack with respect to the $W_{n*}G$ action on $W_{n*}U$ for smooth affine group schemes. More precisely, we prove the following.
	\begin{Thm}
		Let $G\to S$ be a smooth affine group scheme acting on an $S$-scheme $U$, then 
		\begin{center}
			$(J^{n}[U/G]^{pre})^{st}
			\simeq
			[W_{n*}U/W_{n*}G]$
		\end{center}
		
	\end{Thm}
	
	In \cite{article3} Borger showed that the generalized arithmetic jet space functor $W_{n*}$ is the right adjoint of the generalized Witt vector functor $W_{n}^{*}$. In the case of prestacks, category of morphisms between them forms a groupoid. For a scheme $T$ and a stack $\mathscr{X}$, we have the following isomorphism of groupoids

	\begin{Thm}\label{intro_ad}
		Let $T$ be a scheme and $\mathscr{X}$ be a prestack over $(Sch/S)_{\acute{e}t}$, then
		
		\begin{center}
			$Hom(T,J^{n}\mathscr{X})\simeq Hom(W_{n}^{*}T,\mathscr{X})$
		\end{center}
		upto 2-isomorphism of morphisms of stacks over $S$.
	\end{Thm}
	
	In \cite{article3} we have the fact that $W_{n*}$ preserves algebraic spaces. For stacks, we look into algebraic and Deligne-Mumford stacks and we have
	
	\begin{Thm}
		Let $\mathscr{X}$ be an algebraic stack or a Deligne Mumford stack, then $J^{n}\mathscr{X}$ is algebraic stack or a Deligne Mumford stack respectively.
	\end{Thm}
	
	\noindent\textbf{\normalsize Acknowledgments.}
	I would like to thank my advisor Arnab Saha for insightful discussions and Council of Scientific and Industrial Research (CSIR) for the funding under the scheme 09/1031(14019)/2022-EMR-I.
	\section{Preliminaries}
	\label{sec:Preliminaries}
	
	\subsection{Generalized Witt vectors.}	
	We are repeatedly using the notion of the generalized Witt vector sheaves introduced by Borger in \cite[Section 10]{article3} throughout this article. Here we will recall the construction of generalized Witt vector ring, leading to the construction of the generalized Witt vector sheaves. The construction of infinite length generalized Witt vector functor from the category of $R$-algebras $\text{Ring}_{R}$ to itself in two steps. Here we will be discussing an overview of the construction that is given in \cite[Section 1]{article2}.
	\vspace{0.2cm}
	\\
	\textbf{Step 1}: Construction of $W_{R,E}^{\text{fl}}$
	\begin{definition}[Supramaximal ideals]\cite[1.1]{article2}
		An ideal $\mathfrak{m}$ of a ring $R$ is supramaximal if one of the following holds
		\begin{enumerate}
			\item $R/\mathfrak{m}$ is a finite field, $R_{\mathfrak{m}}$ is a discrete valuation ring, and $\mathfrak{m}$ is finitely presented as an $R$-module, or
			\item $\mathfrak{m}$ is a unit ideal.
		\end{enumerate}
		For the rest of the section we are fixing a ring $R$ and a family of pairwise coprime supramaximal ideals of $R$, denoted by $(\mathfrak{m}_{\alpha})_{\alpha \in E}$, where $E$ is an indexing set.\\
		Let $\text{Ring}_{R}^{\text{fl}}$ denote the full subcategory of $\text{Ring}_{R}$, whose objects are $\mathfrak{m}$-torsion free for all $\mathfrak{m}\in E$. We will construct a functor $W_{R,E}^{\text{fl}}:\text{Ring}_{R}^{\text{fl}}\rightarrow \text{Ring}_{R}^{\text{fl}} $. Let $\mathbb{N^{(E)}}:= \bigoplus _{E}\mathbb{N}$ where $\mathbb{N} = \{0,1,2,...\}$ under addition. $\mathbb{N^{(E)}}$ has a commutative monoid structure with coordinate wise addition.
		
		\begin{definition}[$\Psi$-rings]\cite[1.5]{article2}
			An $R$-algebra $A$ is called a $\Psi _{R,E}$ ring if there is an action of the monoid $\mathbb{N} ^{(E)}$.
		\end{definition} 
		
		\begin{exam}\cite[1.5]{article2}
			Consider $R[x]^{\otimes _{R}\mathbb{N} ^{(E)}}$, where the action of $\mathbb{N} ^{(E)}$ is by taking the exponential. It is a free $\Psi_{R,E} $-ring on one generator.
		\end{exam}
		
	\end{definition}
	
	\begin{definition}[$E$-flat $R$-modules.]\cite[1.3]{article2}
		An $R$-module $M$ is $E$-flat if for any maximal ideals $\mathfrak{m}\in E$, we have
		\begin{enumerate}
			\item $R_{\mathfrak{m}}\otimes _{R}M $ is flat over $R_{\mathfrak{m}}$.
			\item the natural map $\mathfrak{m}\otimes _{R}M\rightarrow M$ is injective.
		\end{enumerate}
		
	\end{definition}
	
	\begin{definition}[E-flat $\lambda$-rings ]\cite[1.6]{article2}
		Let $A$ be an $E$-flat $R$-algebra. A $\Lambda _{R,E}$- ring structure on $A$ is a $\Psi_{R,E}$-ring structure, with the frobenius lift property for any $\alpha \in E $ as shown below
		\begin{center}
			$Id\otimes \psi _{\alpha}:R/\mathfrak{m}\otimes A \rightarrow R/\mathfrak{m}\otimes A$
		\end{center} 
		agrees with the Frobenius map $x\mapsto x^{q^{\alpha}}$. Here the morphism between two $\Lambda$-ring is a morphism of $\Psi_{R,E}$-rings. Note that $\psi_{\alpha}$ is the endomorphism of $A$ induced by the action of $\alpha^{th}$ basis vector of $\mathbb{N^{(E)}}$.
	\end{definition}
	
	\begin{definition}[$E$-flat Witt vector rings]\cite[1.8]{article2}
		Consider an $E$-flat $R$ algebra $A$,
		\begin{center}
			$U_{0}(A):= A^{\mathbb{N^{(E)}}}$
		\end{center}
		where $A^{\mathbb{N^{(E)}}}$ is a cofree $\Psi _{R,E}$ ring.	Now inductively we define, for $i \geq 0$
		\begin{center}
			$U_{i+1}(A): = \{ b\in U_{i}(A)| \psi_{\alpha} (b)-b^{q_{\alpha}}\in \mathfrak{m}_{\alpha}U_{i}(A),\text{for all}\enspace \alpha \in E\}$
		\end{center}
		We define the $E$-flat Witt vector ring as
		\begin{center}
			$W^{\text{fl}}_{R,E}(A): = \bigcap_{i\geq 0}U_{i}(A)$
		\end{center}
		which is an $R$-sub algebra of $ A^{\mathbb{N^{(E)}}}$.
	\end{definition}
	\begin{rmk}
		From \cite[Proposition 1.9]{article2}, $W^{\text{fl}}_{R,E}(A)$ have a $\Lambda_{R,E}$-ring structure and the functor $A \mapsto W^{\text{fl}}_{R,E}(A)$ is right adjoint of the forgetful functor.
	\end{rmk}
	\flushleft
	It turns out that the functor $W^{\text{fl}}_{R,E}$ is representable by a flat $R$-algebra $\Lambda_{R,E}$. See \cite[1.11]{article2} for detailed construction. We will have a glimpse at the construction of $\Lambda_{R,E}$. Let $\Lambda_{R,E}^{0} = \Psi_{R,E} $. Now assume that the ring $\Lambda_{R,E} ^{i}$ has been constructed, which is a sub-R-algebra of
	\begin{center}
		$R[1/E]\otimes _{R}\Lambda_{R,E}^{i} = R[1/E]\otimes _{R}\Psi_{R,E}$
	\end{center}
	Now we inductively define
	\begin{center}
		$\Lambda_{R,E} ^{i+1} : =< \{\pi^{*}\otimes(\psi_{\alpha}(f)-f^{q_{\alpha}})| \pi ^{*}\in \mathfrak{m}_{\alpha}^{-1}\subseteq R[1/E],f\in \Lambda_{R,E}^{i},\alpha \in E\}>$
	\end{center} 
	where $<C>$ denote the $sub-\Lambda_{R,E} ^{i}$-algebra of $R[1/E]\otimes _{R}\Psi_{R,E}$ generated by $C\subseteq R[1/E]\otimes _{R}\Psi_{R,E}$.
	We now set,
	\begin{center}
		$\Lambda_{R,E} = \bigcup _{i\geq 0}  \Lambda_{R,E} ^{i}$
	\end{center}
	
	Now since $W^{\text{fl}}_{R,E}$ is representable by a flat $R$- algebra $\Lambda_{R,E}$, we can define a functor extending $W^{\text{fl}}_{R,E}$ as
	\begin{center}
		$W_{R,E}(A) = \text{Hom}_{\text{Ring}_{R}}(\Lambda_{R,E},A)$.
	\end{center}
	resulting in the functor $W_{R,E}: \text{Ring}_{R}\rightarrow \text{Ring}_{R}$.
	
	\subsection{Truncated Witt vectors.}
	Now we give an outline of the truncated Witt vector rings, which are constructed in \cite[Section 2]{article2}.
	\begin{definition}[Ordering]
		Let $n,n' \in \mathbb{Z}^{(E)} = \bigoplus _{E} \mathbb{Z}$, we say $n'\leq n$ if $n'_{\alpha}\leq n_{\alpha},\forall \alpha \in E.$ 
	\end{definition}
	We denote
	\begin{center}
		$[0,n] =\{n'\in \mathbb{N}^{(E)}|\enspace n'\leq n\} $.
	\end{center}
	
	The ring $\Psi_{R,E}$ has the following decomposition,
	\begin{align}\label{decomposition}
		\Psi_{R,E}
		= \bigotimes\limits_{n \in \mathbb{N}^{(E)}} R[\psi_n| n\in \mathbb{N^{(E)}}].
	\end{align}
	Let $n\in \mathbb{Z}^{(E)}$, define
	\begin{align}
		\Psi_{R,E,n} = \bigotimes\limits_{n \in [0,n]} R[\psi_n]
	\end{align}
	which induces the following filtration on $\Lambda_{R,E}$,
	\begin{align}
		\Lambda_{R,E,n} = \Lambda_{R,E}\cap (R[1/E]\otimes _{R}\Psi_{R,E,n})
	\end{align}
	Now we are ready to define the Witt vector rings of finite length as,
	\begin{align}
		W_{R,E,n}(A) = \text{Hom}_{R}(\Lambda_{R,E,n},A)
	\end{align}
	the functor $	W_{R,E,n} $ is called \text{E-typical Witt vector functor of finite length.}

	\subsection{Preliminaries on stacks}
	The theory of stacks has its roots in moduli theory, which serves as a framework for understanding isomorphism classes of families of algebraic geometric objects such as vector bundles, curves, principal bundles etc. in a systematic way. Here we will recall some definitions involving stacks which will be used in the later sections.
	
	\begin{definition}[Prestacks]\cite[Definition~2.4.1]{textbook2}
		A functor $p:\mathscr{X}\rightarrow \mathcal{S}$ is a prestack over a category $\mathcal{S}$ if
		\begin{enumerate}
			\item (pullbacks exits) for every diagram 
			\begin{center}
				$	\xymatrix{a \ar@{.>}[r] \ar@{|.>}[d]& b\ar@{|->}[d]\\ S \ar[r] &T} $
			\end{center}
			
			\item (Universal property for pullbacks) for every diagram
			\begin{center}
				$	\xymatrix{a \ar@{.>}[r] \ar@{|->}[d]\ar@/^/[rr] &b\ar@{->}[r]\ar@{|->}[d] &c\ar@{|->}[d]\\ R \ar@{->}[r]&S \ar[r] &T} $
			\end{center}
			of solid arrows, there exist unique morphism $a\rightarrow b$ over $R\rightarrow S$ filling the diagram.
		\end{enumerate}
	\end{definition}
	
	From now on we write $\mathscr{X}$ for a prestack, rather than giving the data $p:\mathscr{X}\rightarrow \mathcal{S}$. 
	
	\begin{definition}[Stacks]\cite[Definition~2.5.1]{textbook2}
		A prestack $\mathscr{X}$ over a site $\mathcal{S}$ is a stack if the following conditions for all coverings $\{S_{i}\rightarrow S\}$ of an object $S\in \mathcal{S}$:
		\begin{enumerate}
			\item (morphisms glue) For objects $a$ and $b$ in $\mathscr{X}$ over $S$ and morphisms $\phi _{i}:a|_{S_{i}}\rightarrow b$ such that $\phi_{i}|_{S_{ij}}= \phi_{j}|_{S_{ij}}$ as displayed in the diagram 
			\begin{center}
				\begin{tabular}{c c c}
					\xymatrix{
						& a|_{S_i} \ar@/^1.0pc/[rrd]^{\phi_{i}} \ar[rd] \\
						a|_{S_{ij}} \ar[ur] \ar[dr] && a \ar@{.>}[r] & b \\
						& a|_{S_j} \ar@/_1pc/[urr]^{\phi_j} \ar[ur]
					}
					
					&
					
					\raisebox{-1.5cm}{over} 
					
					&
					
					\xymatrix{
						& S_i \ar[dr] & \\
						S_{ij} \ar[ur] \ar[dr] & & S \\
						& S_j \ar[ur] &
					}
				\end{tabular}
				
			\end{center}
			there exist a unique morphism $\phi : a\rightarrow b$ over $Id_{S}$ with $\phi |_{S_{i}} = \phi_{i}$ 
			
			\item (objects glue) For objects $a_{i}$ over $S_{i}$ and isomorphisms $\alpha _{i,j}:a_{i}|_{S_{ij}}\rightarrow a_{j}|_{S_{ij}}$ as displayed in the diagram
			
			\begin{center}
				\begin{tabular}{c c c}
					$\xymatrix@C=1.2em@R=1.8em{
						& a_i \ar[rrdd] & \\
						a_i|_{S_{ij}} \ar[ur] \ar[dd]_{\alpha_{i,j}} & & \\
						& & &a \\
						a_j|_{S_{ij}} \ar[rd] & & \\
						& a_j \ar[uurr]
					}$
					
					&
					
					\raisebox{-2.3cm}{over}
					
					&
					
					\raisebox{-1.1cm}{$
						\xymatrix@C=1.2em@R=1.8em{
							& S_i \ar[dr] & \\
							S_{ij} \ar[ur] \ar[dr] & & S \\
							& S_j \ar[ur] &
						}$}
				\end{tabular}		
			\end{center}
			satisfying the cocycle condition $\alpha _{jk}|_{S_{ijk}}\circ  \alpha _{ij}|_{S_{ijk}}= \alpha _{ik}|_{S_{ijk}}$ on $S_{ijk}$, there exists an object $a$ over $S$ and isomorphisms $\phi _{i}:a|_{S_{i}}\rightarrow a_{i}$ over $Id_{S_{i}}$ such that $\phi_{j}|_{S_{ij}} = \alpha_{i,j}\circ\phi_{i}|_{S_{ij}}$ on $S_{ij}$.
			
		\end{enumerate}
	\end{definition}
	
	Throughout the discussion, an algebraic space is,
	\begin{definition}\cite[Definition~5.1.10]{textbook1}
		Let $S$ be a scheme. An algebraic space over $S$ is a functor $X:(Sch/S)^{op}\rightarrow Set $ such that the following hold:
		\begin{enumerate}
			\item X is a sheaf with respect to the big étale topology.
			\item $\Delta : X\rightarrow X\times _{S} X$ is representable by schemes.
			\item There exists an $S$-scheme $U\rightarrow S$ and a surjective étale morphism $U\rightarrow X$
		\end{enumerate}
		
	\end{definition}

	\begin{definition}\cite[Definition~3.1.3]{textbook2}
		A morphism $\mathscr{X}\rightarrow \mathscr{Y}$ of prestacks over Sch is representable if for every morphism $T\rightarrow \mathscr{Y}$ from a scheme $T$, the fiber product $\mathscr{X}\times _{\mathscr{Y}} T$ is an algebraic space.
	\end{definition}

	\begin{definition}\cite[Definition~8.1.4]{textbook1}
		A stack $\mathscr{X}/S$ is an algebraic stack if the following hold:
		\begin{enumerate}
			\item The diagonal 
			\begin{center}
				$\Delta :\mathscr{X}\rightarrow \mathscr{X}\times _{S}\mathscr{X}$
			\end{center}
			is representable.
			
			\item There exist a smooth surjective morphism $\pi :X\rightarrow \mathscr{X}$ with $X$ a scheme.
		\end{enumerate}
	\end{definition}
	
	\begin{definition}\cite[Definition~8.3.1]{textbook1}
		
		An algebraic stack $\mathscr{X}/S$ is called Deligne- Mumford if there exists an étale surjection $X\rightarrow \mathscr{X}$ with $X$ a scheme.
	\end{definition}
	
	\section{Construction of jet spaces on stacks.}
	
	\label{sec: COnstruction of jets}
	Consider the category $(Sch/S)_{\acute{e}t}$ of schemes over a base scheme $S$ with  étale topology equipped on it. We have the following assumptions on $S$.
	\begin{enumerate} [label=(\roman*)]
		\item When $S$ is affine, then $S$ is the spectrum of a ring $R$ with properties specified in \cite[~1.1]{article2} or Section \ref{sec:Preliminaries}.
		\item \cite[10.2]{article3} In general, $S$ is assumed to be a separated scheme. A supra-maximal ideal on $S$ to be a finitely presented ideal sheaf $\mathfrak{m}$ in $\mathcal{O} _{S}$ corresponding to either\\
		(a) a closed point whose local ring $\mathcal{O}_{S,\mathfrak{m}}$ is a discrete valuation ring with finite residue filed, or \\
		(b) the empty subscheme.
	\end{enumerate}
	Throughout this paper $E$ is a sub-collection of supra-maximal ideals of $S$, which are pairwise coprime and $S,E,n$ is arbitrarily fixed. The generalized Witt vector functor $W^{*}_{E,S,n}$ defined on sheaves over $S$ is considered for defining the jet space functor on prestacks. From now on we will not specify $E,S$ and will denote $W^{*}_{E,S,n}$ as $W^{*}_{n}$.\\
	\textbf{Step 1}: (Construction of $\mathscr{W}_{n}^{*}$)\\
	We define a category $\mathscr{W}_{n}^{*}$ whose objects are given by $W_{n}^{*}U$ for a given $S$-scheme $U$.  A morphism $W_{n}^{*}V \rightarrow W_{n}^{*}U $ exists if there is a morphsim of schemes $f: V\rightarrow U $, and the unique morphism is given by $W_{n}^{*}(f)$. \\
	From \cite[Corollary 15.8]{article3}, we know  that the functor $W_{n}^{*}$ is faithful from the category of algebraic spaces to itself and $W_{n}^{*}$ maps schemes to schemes by \cite[Corollary 15.6]{article3}. Since schemes are algebraic spaces, we have the following faithfull functor between categories 
	
	\begin{align}
		W_{n}^{*}: (Sch/S)\rightarrow (Sch/S)
	\end{align}
	By construction of $\mathscr{W}_{n}^{*}$, we have 
	
	\begin{center}
		$Hom_{(Sch/S)}(V,U)\simeq Hom_{\mathscr{W}_{n}^{*}}(W_{n}^{*}V,W_{n}^{*}U)$
	\end{center}
	$W_{n}^{*}: (Sch/S) \rightarrow \mathscr{W^{*}}_{n}$ is the restriction to the image of the functor defined in (1.1) and by construction it is fully faithful and essentially surjective, giving rise to an equivalence of categories
	
	\begin{center}
		$(Sch/S) \simeq \mathscr{W^{*}}_{n}$ 
	\end{center}
	So, we have a functor $p' : \mathscr{W}_{n}^{*}\rightarrow (Sch/S)$ such that $W_{n}^{*}U\mapsto U$ and in the case of morphisms, ($W_{n}^{*}(f):W_{n}^{*}V \rightarrow W_{n}^{*}U)\mapsto (f: V\rightarrow U) $.
	\vspace{0.2cm} \\
	\textbf{Step 2} : (Construction of jet prestack)
	\vspace{0.2cm}
	\\
	Consider a prestack $(\mathscr{X},p)$ over the category $(Sch/S)$. \\
	We define a category $J^{n}\mathscr{X}$, whose objects are defined as follows, given an $S$-scheme $U$ the objects are given by $a\in \mathscr{X}(W_{n}^{*}U)$.  Given two objects $ a,b \in Ob(J^{n}\mathscr{X})$, a morphism $\alpha : a\rightarrow b $ in $\mathscr{X}$ is a morphism in $J^{n}\mathscr{X}$ if there is $f:V\rightarrow U$ such that $p(\alpha) = W_{n}^{*}(f):p(a) =W_{n}^{*}V \rightarrow p(b)= W_{n}^{*}U $. \\
	Now we have to define a functor $\tilde{p}:J^{n}\mathscr{X}\rightarrow (Sch/S)$, which gives the prestack structure.\\   
	\textbf{Step 3}: (Constructing the functor)
	\vspace{0.2cm}
	\\
	There is a natural restriction functor induced from $p$, by abuse of notation we call it $p$ itself, which is given by 
	
	\begin{center}
		$p:J^{n}\mathscr{X}\rightarrow \mathscr{W}_{n}^{*}$ 
	\end{center}
	such that for $a\in Ob(J^{n}\mathscr{X})$, that is $a\in \mathscr{X}(W_{n}^{*}(U))$ for some scheme $U$ then $p(a) = W_{n}^{*}U$. Now given a morphism  $\alpha :a \rightarrow b $ in  $J^{n}\mathscr{X}$, by the definition $p(\alpha) =  J^{n}(f): J^{n}V\rightarrow J^{n}U$ for some $f: V\rightarrow U $, which is a genuine morphism in $ \mathscr{W}_{n}^{*}$.
	\newline
	We have defined two categories namely $J^{n}\mathscr{X}$, $\mathscr{W}_{n}^{*}$ and functors $p:J^{n}\mathscr{X}\rightarrow \mathscr{W}_{n}^{*}$ and $p' : \mathscr{W}_{n}^{*} \rightarrow (Sch/S)$. Consider the functor $\tilde{p}:J^{n}\mathscr{X} \rightarrow (Sch/S)$, where $\tilde{p} = p'\circ p$. Finally, the data $(J^{n}\mathscr{X},\tilde{p})$ is the candidate for extending the arithmetic jet space functor to stacks. 
	\begin{Lem}
		$(\mathscr{W}_{n}^{*},p')$ is a prestack over the category $(Sch/S)$.	
	\end{Lem}
	\begin{proof}
		(Pull back exists)\\
		Give a diagram of solid arrows as below, we have to find an object and a morphism making the square complete.
		\begin{center}
			$	\xymatrix{ &W_{n}^{*}U\ar@{|->}[d]^{p'}\\ V \ar[r]^f &U} $
		\end{center}
		Consider the object $W_{n}^{*}V$, the unique morphism $W_{n}^{*}(f)$ in the category $\mathscr{W}_{n}^{*}$ completes the above square and we have,
		
		\begin{center}
			$	\xymatrix{W_{n}^{*}V\ar@{.>}[r]^{W_{n}^{*}(f)}\ar@{|.>}[d]^{p'} &W_{n}^{*}U\ar@{|->}[d]^{p'}\\ V \ar[r]^f &U} $
		\end{center}
		(Universal property of pullbacks)\\
		Given  a diagram of solid arrows as below, 
		
		\begin{center}
			$	\xymatrix{W_{n}^{*}R \ar@{|->}[d]^{p'}\ar@/^1pc/[rr]^{W_{n}^{*}(f\circ g)} &W_{n}^{*}V\ar@{->}[r]_{W_{n}^{*}(f)}\ar@{->}[d]^{p'} &W_{n}^{*}U\ar@{|->}[d]^{p'}\\ R \ar@{->}[r]^{g}&V \ar[r]^f &U} $
		\end{center}
		There is a unique morphism between $W_{n}^{*}R$ and $W_{n}^{*}V$ making the above diagram functorial. Uniqueness is immediate from the construction of the category $\mathscr{W}_{n}^{*}$. Note that $W_{n}^{*}$ is functorial, hence $W_{n}^{*}(g):W_{n}^{*}R \rightarrow W_{n}^{*}V$ fits into the diagram since $W_{n}^{*}(f\circ g) = W_{n}^{*}(f)\circ W_{n}^{*}(g)$. We have the completed diagram as,
		
		\begin{center}
			$	\xymatrix{W_{n}^{*}R\ar@{.>}[r]_{W_{n}^{*}(g)} \ar@{|->}[d]^{p'}\ar@/^1pc/[rr]^{W_{n}^{*}(f\circ g)} &W_{n}^{*}V\ar@{->}[r]_{W_{n}^{*}(f)}\ar@{|->}[d]^{p'} &W_{n}^{*}U\ar@{|->}[d]^{p'}\\ R \ar@{->}[r]^{g}&V \ar[r]^f &U} $
		\end{center}	
	\end{proof}
	
	\begin{Prpn}
		Given a prestack $(\mathscr{X},p)$ then $(J^{n}\mathscr{X},\tilde{p})$ is a prestack over  $Sch/S$.
	\end{Prpn}
	
	\begin{proof}

		Claim 1: (Pull back  exists)
		
		\begin{center}
			$	\xymatrix{ & a\ar@{|->}[d]^{\tilde{p}}\\ V \ar[r]^f &U} $
		\end{center}
		As we have done before, we have to complete the above functorial diagram by finding an object and a morphism in the category $J^{n}\mathscr{X}$ fitting the above data. Unwinding the above diagram we have,
		\begin{center}
			$	\xymatrix{& a \ar@{|->}[d]^{p}\\W_{n}^{*}V\ar@{.>}[r]^{W_{n}^{*}(f)}\ar@{|.>}[d]^{p'} &W_{n}^{*}U\ar@{|->}[d]^{p'}\\ V \ar[r]^f &U} $
		\end{center}
		
		now using the fact that $(\mathscr{X},p)$ is a prestack, $\exists b\in Ob(J^{n}\mathscr{X})$ such that $p(x) = W_{n}^{*}V $. We have the completed diagram as, 
		
		\begin{center}
			$	\xymatrix{b \ar@{|.>}[d]^{p} \ar@{.>}[r]  & a \ar@{|->}[d]^{p}\\W_{n}^{*}V\ar@{.>}[r]^{W_{n}^{*}(f)}\ar@{|.>}[d]^{p'} &W_{n}^{*}U\ar@{|->}[d]^{p'}\\ V \ar[r]^f &U} $
		\end{center}
		
		Claim 2 : (Universal property of pullbacks)\\
		\vspace{0.3cm}
		Consider the diagram given below,

		\begin{center}
			$	\xymatrix{c \ar@{|->}[d]^{\tilde{p}}\ar@/^1pc/[rr] &b\ar@{->}[r]\ar@{|->}[d]^{\tilde{p}} &a\ar@{|->}[d]^{\tilde{p}}\\ R \ar@{->}[r]^{g}&V \ar[r]^f &U} $
		\end{center}
		We have to show that there is a unique morphism from $c$ to $b$ making the diagram functorial. Again unwinding the diagram using the definition of $\tilde{p}$, 
		
		\begin{center}
			$	\xymatrix{c \ar@{|->}[rd]^{p}\ar@/^1pc/[rr] &b\ar@{->}[r]\ar@{|->}[rd]^{p} &a\ar@{|->}[rd]^{p}\\&W_{n}^{*}R \ar@{|->}[ld]_{p'}\ar@/^1pc/[rr]^{W_{n}^{*}(f\circ g)} &W_{n}^{*}V\ar@{->}[r]_{W_{n}^{*}(f)}\ar@{|->}[ld]_{p'} &W_{n}^{*}U\ar@{|->}[ld]^{p'}\\ R \ar@{->}[r]^{g}&V \ar[r]^f &U} $
		\end{center}
	    Using the fact that $\mathscr{W}_{n}^{*}$ is a prestack, there is a unique morphism $W_{n}^{*}(g)$ making the middle diagram commute. Now by using the fact that $\mathscr{X}$ is a prestack there is a unique morphism from $c$ to $b$ making the upper part of the diagram commute. Now we have the completed diagram as,
		\begin{center}
			$	\xymatrix{c \ar@{.>}[r]\ar@{|->}[rd]^{p}\ar@/^1pc/[rr] &b\ar@{->}[r]\ar@{|->}[rd]^{p} &a\ar@{|->}[rd]^{p}\\&W_{n}^{*}R \ar@{.>}[r]_{W_{n}^{*}(g)}\ar@{|->}[ld]_{p'}\ar@/^1pc/[rr]^{W_{n}^{*}(f\circ g)} &W_{n}^{*}V\ar@{->}[r]_{W_{n}^{*}(f)}\ar@{|->}[ld]^{p'} &W_{n}^{*}U\ar@{|->}[ld]^{p'}\\ R \ar@{->}[r]_{g}&V \ar[r]_f &U} $
		\end{center}

	\end{proof}

	\begin{Thm}
		If $\mathscr{X}$ is a stack then the prestack $J^{n}\mathscr{X}$ is a stack.
	\end{Thm}
	\begin{proof}
		Consider the prestack $J^{n}\mathscr{X}$, let $U\in (Sch/S)$ and $\{U_{i}\rightarrow U\}$ be an étale covering. \\
		(1) Morphisms glue. \\
		Let $a,b\in Ob(J^{n}\mathscr{X})$ and $\phi _{i} : a|_{U_{i}}\rightarrow b $ be such that $\phi _{i}|_{U_{ij}} =\phi _{j}|_{U_{ij}}$, which is shown in the below diagram 
		
		\begin{center}
			\begin{tabular}{c c c}
				\xymatrix{
					& a|_{U_i} \ar@/^1.0pc/[rrd]^{\phi_{i}} \ar[rd] \\
					a|_{U_{ij}} \ar[ur] \ar[dr] && a \ar@{.>}[r] & b \\
					& a|_{U_j} \ar@/_1pc/[urr]^{\phi_j} \ar[ur]
				}
				
				&
				
				\raisebox{-1.5cm}{over} 
				
				&
				
				\xymatrix{
					& U_i \ar[dr] & \\
					U_{ij} \ar[ur] \ar[dr] & & U \\
					& U_j \ar[ur] &
				}
			\end{tabular}
			
		\end{center}
		
		we will show the existence of a unique morphism $\phi : a\rightarrow b $ over $id_{S}$ with $\phi |_{S_{i}} = \phi _{i}$.
		By the definition of $J^{n}\mathscr{X}$ the diagram above can be visualized as,
		\begin{center}
			\begin{tabular}{c c c}
				\xymatrix{
					& a|_{W_{n}^{*}U_i} \ar@/^1.0pc/[rrd]^{\phi_{i}} \ar[rd] \\
					a|_{W_{n}^{*}U_{ij}} \ar[ur] \ar[dr] && a \ar@{.>}[r] & b \\
					& a|_{W_{n}^{*}U_j} \ar@/_1pc/[urr]^{\phi_j} \ar[ur]
				}
				
				&
				
				\raisebox{-1.5cm}{over} 
				
				&
				
				\xymatrix{
					& W_{n}^{*}U_i \ar[dr] & \\
					W_{n}^{*}U_{ij} \ar[ur] \ar[dr] & & W_{n}^{*}U \\
					& W_{n}^{*}U_j \ar[ur] &
				}
			\end{tabular}
			
		\end{center}
		Note that this diagram is now in the stack $\mathscr{X}$. By \cite[Corollary 15.3]{article3} for an étale cover $(U_{i})_{i\in I}$ of a scheme $U$ over $S$ then $(W_{n}^{*}U_{i})_{i\in I}$ is an étale cover of $W_{n}^{*}X$ and 
		\begin{center}
			$W_{n}^{*}(U_{i}\times _{U}U_{j})\rightarrow W_{n}^{*}(U_{i})\times _{W_{n}^{*}(U)}W_{n}^{*}(U_{j})$ given by $(W_{n}^{*}(pr_{1}),W_{n}^{*}(pr_{2}))$  
		\end{center} 
		is an isomorphism. We denote by $\tilde{U}_{ij}$ for $W_{n}^{*}(U_{i})\times _{W_{n}^{*}(U)}W_{n}^{*}(U_{j})$.
		\begin{align}\label{glueing1}
			\xymatrix{
				&& a|_{W_{n}^{*}U_i} \ar@/^1.0pc/[rrd]^{\phi_{i}} \ar[rd] \\
				a|_{\tilde{U}_{ij}}\ar@{.>}[r] \ar@{.>}[urr]\ar@{.>}[drr]&a|_{W_{n}^{*}U_{ij}} \ar[ur] \ar[dr] && a \ar@{.>}[r] & b \\
				&& a|_{W_{n}^{*}U_j} \ar@/_1pc/[urr]^{\phi_j} \ar[ur]	}		
		\end{align}	
		In the diagram above, the morphism
		\[
		a|_{\widetilde{U}_{ij}} \longrightarrow a|_{W_n^{*}U_{ij}}
		\]
		is obtained by pullback along the isomorphism
		\[
		\bigl(W_n^{*}(\mathrm{pr}_1),\, W_n^{*}(\mathrm{pr}_2)\bigr).
		\]
		The morphisms
		\[
		a|_{\widetilde{U}_{ij}} \longrightarrow a|_{W_n^{*}U_i}
		\quad\text{and}\quad
		a|_{\widetilde{U}_{ij}} \longrightarrow a|_{W_n^{*}U_j}
		\]
		are the unique morphisms obtained by invoking the axioms of prestacks. From the diagram, we obtain
		\[
		\phi_i|_{\widetilde{U}_{ij}} = \phi_j|_{\widetilde{U}_{ij}},
		\]
		i.e.\ by commutativity of the diagram, the morphism
		\(\phi_i|_{\widetilde{U}_{ij}}\) is obtained by composing
		\(\phi_i|_{W_n^{*}U_{ij}}\) with the arrow
		\[
		a|_{\widetilde{U}_{ij}} \longrightarrow a|_{W_n^{*}U_{ij}}.
		\]
		\flushleft	
		the same argument is true for 	\(\phi_j|_{\widetilde{U}_{ij}}\). Now, since \(\mathscr{X}\) is a stack, we can invoke the gluing axiom for morphisms, hence there exists a unique morphism
		\[
		\phi \colon a \longrightarrow b
		\]
		over \(\mathrm{id}_{W^n(U)}\) such that
		\(\phi|_{U_i} = \phi_i\) for all \(i\). Noting that
		\[
		p'\bigl(\mathrm{id}_{W^n(U)}\bigr) = \mathrm{id}_U,
		\]
		it follows that \(\phi\) is a genuine morphism in the prestack
		\(J^n\mathscr{X}\).
		\\
		(2) (Objects glue) \\
		Let $a_{i}$ be over $U_{i}$ and isomorphisms $\alpha _{ij}:a_{i}|_{U_{ij}}\rightarrow a_{j}|_{U_{ij}}$ displayed in the diagram,
		\begin{center}
			\begin{tabular}{c c c}
				$\xymatrix@C=1.2em@R=1.8em{
					& a_i \ar[rrdd] & \\
					a_i|_{U_{ij}} \ar[ur] \ar[dd]_{\alpha_{ij}} & & \\
					& & &a \\
					a_j|_{U_{ij}} \ar[rd] & & \\
					& a_j \ar[uurr]
				}$
				
				&
				
				\raisebox{-2.3cm}{over}
				
				&
				
				\raisebox{-1.1cm}{$
					\xymatrix@C=1.2em@R=1.8em{
						& U_i \ar[dr] & \\
						U_{ij} \ar[ur] \ar[dr] & & U \\
						& U_j \ar[ur] &
					}$}
			\end{tabular}		
		\end{center}

		satisfying the cocycle condition 
		\begin{center}
			$\alpha _{jk}|_{U_{ijk}}\circ \alpha _{ij}|_{U_{ijk}} = \alpha _{ik}|_{U_{ijk}}$
		\end{center}
		
		we will show the existence of $a$ over $U$ and isomorphism $\phi _{i}: a|_{U_{i}}\rightarrow a_{i}$ over $Id_{U_{i}}$ such that 
		
		\begin{center}
			$\phi _{j}|_{U_{ij}} = \alpha _{ij}\circ \phi _{i}|_{U_{ij}}$
		\end{center}   
		on $U_{ij}$.\\
		Note that for any object $a$ over $U$, the object $a|_{U_{i}}$ over $U_i$ in the category $J^{n}\mathscr{X}$ is same as $a|_{W_{n}^{*}U_{i}}$ in the category $\mathscr{X}$. We have a sequence of morphisms
		
		\begin{align}
			a_i|_{W_n^* U_{ijk}} \to a_i|_{W_n^* U_{ij}} \xrightarrow{\alpha_{ij}} a_j|_{W_n^* U_{ij}}
		\end{align}
		
		in the category $\mathscr{X}$ obtained from 
		
		\begin{align}
			a_i|_{U_{ijk}} \to a_i|_{U_{ij}} \xrightarrow{\alpha_{ij}} a_j|_{U_{ij}}
		\end{align}
		
		in the category $J^{n}\mathscr{X}$.\\
		Using \cite[Corollary.~15.3]{article3} 
		
		\begin{align}
			W_{n}^{*}(U_{i}\times _{U}U_{j})\rightarrow W_{n}^{*}(U_{i})\times _{W_{n}^{*}(U)}W_{n}^{*}(U_{j}) 
		\end{align}
		given by $(W_{n}^{*}(pr_{1}),W_{n}^{*}(pr_{2})) $ is an isomorphism.\\
		For the rest of the proof, we will use $\tilde{U}_{ijk}$ for $ W_{n}^{*}(U_{i})\times _{W_{n}^{*}(U)}W_{n}^{*}(U_{j})\times _{W_{n}^{*}(U)}W_{n}^{*}(U_{k})$ with indices changed accordingly. 
		\[
		\xymatrix{
			a_i|_{W_n^* U_{ijk}} \ar[r]  & 	a_i|_{W_n^* U_{ij}} \ar[r]^{\alpha_{ij}}   & 	a_j|_{W_n^* U_{ij}}  \\
			&a_{i}|_{\tilde{U}_{ij}} \ar@{.>}[u]^{\beta _{i}}        & a_{j}|_{\tilde{U}_{ij}} \ar@{.>}[u]^{\beta _{j}}
		}
		\]
		
		Here $ a_{i}|_{\tilde{U}_{ij}}$ and $a_{j}|_{\tilde{U}_{ij}}$ are obtained from axiom(1) of prestacks by using the isomorphism $(W_{n}^{*}(pr_{1}),W_{n}^{*}(pr_{2})) $   above. From axiom(2) for prestacks, we have a unique arrow completing the square in the above diagram, resulting in 
		
		\[
		\xymatrix{
			a_i|_{W_n^* U_{ijk}} \ar[r]  & 	a_i|_{W_n^* U_{ij}} \ar[r]^{\alpha_{ij}}   & 	a_j|_{W_n^* U_{ij}}  \\
			&a_{i}|_{\tilde{U}_{ij}}\ar@{.>}[r]^{\tilde{\alpha}_{ij}} \ar@{.>}[u]^{\beta _{i}}        & a_{j}|_{\tilde{U}_{ij}} \ar@{.>}[u]^{\beta _{j}}
		}
		\]
		
		Using the natural isomorphism $W_{n}^{*}U_{ijk}\rightarrow \tilde{U} _{ijk}$, we have a pull back $a_{i}|_{\tilde{U} _{ijk}}$ along the isomorphism, giving the diagram  
		
		\[
		\xymatrix{
			a_i|_{W_n^* U_{ijk}} \ar[r]  & 	a_i|_{W_n^* U_{ij}} \ar[r]^{\alpha_{ij}}   & 	a_j|_{W_n^* U_{ij}}  \\  a_{i}|_{\tilde{U} _{ijk}} \ar@{.>}[u]^{\gamma _{i}}
			&a_{i}|_{\tilde{U}_{ij}}\ar@{.>}[r]^{\tilde{\alpha}_{ij}} \ar@{.>}[u]^{\beta _{i}}        & a_{j}|_{\tilde{U}_{ij}} \ar@{.>}[u]^{\beta _{j}}
		}
		\]
		
		by invoking axiom(2) of prestacks, we complete the left hand square making the whole diagram commute
		
		\[
		\xymatrix{
			a_i|_{W_n^* U_{ijk}} \ar[r]  & 	a_i|_{W_n^* U_{ij}} \ar[r]^{\alpha_{ij}}   & 	a_j|_{W_n^* U_{ij}}  \\  a_{i}|_{\tilde{U} _{ijk}} \ar@{.>}[u]^{\gamma _{i}} \ar@{.>}[r]
			&a_{i}|_{\tilde{U}_{ij}}\ar@{.>}[r]^{\tilde{\alpha}_{ij}} \ar@{.>}[u]^{\beta _{i}}        & a_{j}|_{\tilde{U}_{ij}} \ar@{.>}[u]^{\beta _{j}}
		}
		\]
		
		Note that in the above diagram, all the horizontal arrows are unique once a pull back is chosen, which are denoted by the dotted vertical arrows. By the same technique used above, we have a commutative diagram 
		
		\[
		\xymatrix{
			a_i|_{W_n^* U_{ijk}} 
			\ar[r] 
			\ar@{.>}@/^1.5pc/[rrr]^{\alpha_{ij}|_{W_n^* U_{ijk}}}
			& a_i|_{W_n^* U_{ij}} \ar[r]^{\alpha_{ij}} 
			& a_j|_{W_n^* U_{ij}} 
			& a_j|_{W_n^* U_{ijk}} \ar@{.>}[l] \\
			a_i|_{\tilde{U}_{ijk}} \ar@{.>}[u]^{\gamma_i} 
			\ar@{.>}[r] 
			\ar@{.>}@/_1.5pc/[rrr]_{\tilde{\alpha}_{ij}|_{\tilde{U}_{ijk}}}
			& a_i|_{\tilde{U}_{ij}} \ar@{.>}[r]^{\tilde{\alpha}_{ij}} \ar@{.>}[u]^{\beta_i} 
			& a_j|_{\tilde{U}_{ij}} \ar@{.>}[u]^{\beta_j} 
			& a_j|_{\tilde{U}_{ijk}} \ar@{.>}[l] \ar@{.>}[u]^{\gamma_j}
		}
		\]
		
		from the above diagram, we have \\
		
		\begin{center}
			$\tilde{\alpha}_{ij}|_{\tilde{U}_{ijk}} = \gamma _{j}^{-1}\circ \alpha_{ij}|_{W_n^* U_{ijk}}\circ \gamma_{i}$.
		\end{center}
		
		following the same procedure for other indices, 
		
		\begin{center}
			$\tilde{\alpha}_{jk}|_{\tilde{U}_{ijk}} = \gamma _{k}^{-1}\circ \alpha_{jk}|_{W_n^* U_{ijk}}\circ \gamma_{j}$.
		\end{center}
		
		\begin{center}
			$\tilde{\alpha}_{ik}|_{\tilde{U}_{ijk}} = \gamma _{k}^{-1}\circ \alpha_{ik}|_{W_n^* U_{ijk}}\circ \gamma_{i}$.
		\end{center}
		Now,\\
		
		\begin{align}
			\tilde{\alpha}_{jk}|_{\tilde{U}_{ijk}} \circ \tilde{\alpha}_{ij}|_{\tilde{U}_{ijk}} 
			&= \gamma_k^{-1} \circ \alpha_{jk}|_{W_n^* U_{ijk}} \circ \underbrace{\gamma_j \circ \gamma_j^{-1}}_{= \mathrm{id}} \circ \alpha_{ij}|_{W_n^* U_{ijk}} \circ \gamma_i \\
			&= \gamma_k^{-1} \circ \alpha_{jk}|_{W_n^* U_{ijk}} \circ \alpha_{ij}|_{W_n^* U_{ijk}} \circ \gamma_i
		\end{align}
		
		but,
		
		\begin{align}
			\alpha_{jk}|_{W_n^* U_{ijk}} \circ \alpha_{ij}|_{W_n^* U_{ijk}}& = \alpha_{ik}|_{W_n^* U_{ijk}}
		\end{align}
		so we have,
		
		\begin{align}
			\tilde{\alpha}_{jk}|_{\tilde{U}_{ijk}} \circ \tilde{\alpha}_{ij}|_{\tilde{U}_{ijk}} 
			&=\tilde{\alpha}_{ik}|_{\tilde{U}_{ijk}}
		\end{align}
		
		since $\mathscr{X}$ is a stack, there exists $a\in \mathscr{X}(W_{n}^{*})$ and isomorphisms $\phi _{i}: a|_{W_{n}^{*}U_{i}}\rightarrow a_{i}$ over $Id _{W_{n}^{*}U_{i}}$, but by the construction, for all $i$, $\phi_{i}$ is a morphism in $J^{n}\mathscr{X}$ over $Id_{U_{i}}$. \\
		Now consider the diagram
		
		\[
		\xymatrix{
			a|_{\tilde{U}_{ij}}\ar@{.>}[d]^{\Theta} \ar[r]^{\phi_{i}|_{\tilde{U}_{ij}}} \ar@{->}@/^2.5pc/[rr]^{\phi _{j}|_{\tilde{U}_{ij}}} & 	a_i|_{\tilde{U}_{ij}}\ar@{.>}[d]^{\beta _{i}} \ar[r]^{\tilde{\alpha}_{ij}}   & 	a_j|_{\tilde{U}_{ij}} \ar@{.>}[d]^{\beta _{j}} \\  a|_{W_{n}^{*}U_{ij}}  \ar@{.>}[r]^{\phi_{i}|_{W_{n}^{*}U_{ij}}}\ar@{.>}@/_2.5pc/[rr]^{\phi _{j}|_{W_{n}^{*}{U}_{ij}}}
			&a_{i}|_{W_{n}^{*}U_{ij}}\ar[r]^{\alpha_{ij}}         & a_{j}|_{W_{n}^{*}U_{ij}} 
		}
		\]
		similar to the previous case, the vertical arrows indicate pull-backs and the dotted horizontal arrow indicates unique arrow making the whole diagram commute.\\
		From the properties of the stack $\mathscr{X}$, we have
		\begin{equation}
			\phi _{j}|_{\tilde{U}_{ij}} = \tilde{\alpha}_{ij}\circ \phi _{i}|_{\tilde{U}_{ij}}
		\end{equation}
		
		Note that the horizontal sequence in the bottom is obtained from the pullback of the vertical arrows and they are commutative, so we have
		
		\begin{equation}
			\phi _{j}|_{W_{n}^{*}{U}_{ij}} = \tilde{\alpha}_{ij}\circ \phi _{i}|_{W_{n}^{*}{U}_{ij}}
		\end{equation}
		
		$\phi _{i}|_{W_{n}^{*}{U}_{ij}}$ and $\phi _{j}|_{W_{n}^{*}{U}_{ij}}$ are morphisms over $id_{{W_{n}^{*}{U}_{ij}}}$, which then becomes genuine morphisms in the category $J^{n}\mathscr{X}$, hence we have an element $a\in J^{n}\mathscr{X}(U)$ satisfying the required properties.
			
	\end{proof}
	Given a stack $\mathscr{X}$, we have constructed a stack $J^{n}\mathscr{X}$. Now we have to establish that its restriction to sheaves agrees with the original construction of generalized arithmetic jet space on sheaves.
	
	\begin{Lem}\label{Sheaf}
		Let $F$ be a sheaf,  $\mathscr{X}_{W_{n*}F}$ be the stack representing the sheaf $W_{n*}F$ and $J^{n}\mathscr{X}_{F}$ be the stack obtained after applying the $J^{n}$ functor to the stack representing the sheaf $F$, then
		\begin{center}
			$\mathscr{X}_{W_{n*}F} \simeq J^{n}\mathscr{X}_{F}$
		\end{center}
		as stacks.
	\end{Lem}
	
	\begin{proof}
		Let $U$ be a scheme over $S$, an object in the stack over $U$ in the stack $ \mathscr{X}_{W_{n*}F}$ is of the form $(a,U)$ where $a \in W_{n*}F(U)$, similarly an object over $U$ in the stack $J^{n}\mathscr{X}_{F}$ is of the form $(b,W_{n}^{*}U)$ where $b\in F(W_{n}^{*}U)$\\
		Now we define a functor, which turns out to be an isomorphism of stacks by invoking the adjunction of functors $W_{n*}$ and $W_{n}^{*}$. Let $f:\mathscr{X}_{W_{n*}F} \rightarrow J^{n}\mathscr{X}_{F}$ be such that
		
		\begin{center}
			$f((a,U)):= (\tilde{a},W_{n}^{*}U)$
		\end{center}
		where $\tilde{a}$ is the adjoint image of $a$ in $X(W_{n}^{*}U)$. Now consider a morphism $(b,V)\xrightarrow{\alpha}(a,U)$, where $\alpha : V\rightarrow U$ and $b = \alpha ^{*}a $. Now we can define
		\begin{center}
			$f((b,V)\xrightarrow{\alpha}(a,U)):=(\tilde{b},W_{n}^{*}V)\xrightarrow{W_{n}^{*}(\alpha)}(\tilde{a},W_{n}^{*}U)$
		\end{center}
		which is an equivalence of stacks by successively using adjunction of the functors $W_{n*}$ and $W_{n}^{*}$.
	\end{proof}
	As the generalized arithmetic jet spaces on sheaves were defined as functor of points using generalized Witt vector rings, the jet functor naturally distribute over the fibre product of sheaves. For stacks we have the following   
	\begin{Prpn}
		Let $\mathscr{X} \xrightarrow{f} \mathscr{Z}$ and $\mathscr{Y}\xrightarrow{g} \mathscr{Z}$ be morphisms of prestacks over $(Sch/S)_{\acute{e}t}$ then
		\begin{center}
			$J^{n}(\mathscr{X}\times _{f,\mathscr{Z},g}\mathscr{Y})\simeq J^{n}\mathscr{X}\times _{J^{n}(f),J^{n}\mathscr{Z},J^{n}(g)}J^{n}\mathscr{Y} $
		\end{center}
	\end{Prpn}
	
	\begin{proof}
		
		Evaluating the prestacks above for a given $S$-scheme $U$, we have
		
		\begin{align}
			J^{n}(\mathscr{X} \times_{f,\mathscr{Z},g} \mathscr{Y})(U) 
			&= (\mathscr{X} \times_{f,\mathscr{Z},g} \mathscr{Y})(W_n^* U) \\
			&= \mathscr{X}(W_n^* U) \times_{f,\mathscr{Z}(W_n^* U),g} \mathscr{Y}(W_n^* U)
			\label{fibreproduct1}
		\end{align}
		
		\begin{align}
			\left(J^{n} \mathscr{X} \times_{J^{n}(f),J^{n}\mathscr{Z},J^{n}(g)} J^{n} \mathscr{Y} \right)(U) 
			&= J^{n} \mathscr{X}(U) \times_{J^{n}(f),J^{n}\mathscr{Z},J^{n}(g)} J^{n} \mathscr{Y}(U)\\
			&= \mathscr{X}(W_n^*(U)) \times_{f,\mathscr{Z}(W_n^*(U)),g} \mathscr{Y}(W_n^*(U))
			\label{fibreproduct2}
		\end{align}
		from \ref{fibreproduct1} and \ref{fibreproduct2} one observes that by construction, the given prestacks are not just isomorphic, they are same as categories.
	\end{proof}
	
	\section{Jet spaces and quotient stacks}
	 We will consider schemes with a group action on it, the given action induces a natural action in the level of jet spaces via functor of points. Our main objective is to look at the arithmetic jet of the quotient stack obtained by an action of a smooth affine group scheme on a given scheme.

	\begin{definition}[Quotient prestacks]\cite[Definition.~2.4.14]{textbook2}
		
		Let $G\rightarrow S$ be a smooth affine group scheme acting on a scheme $U$ over $S$. The quotient prestacks $[U/G]^{pre}$ of an action of a smooth affine group scheme $G\rightarrow S$ on an $S$-scheme $U$ is the category over $Sch/S$ consisting of pairs $(T,u)$ where $T$ is an $S$-scheme and $u\in U(T)$. A morphism $(T',u')\rightarrow (T,u)$ is the data of a map $f: T'\rightarrow T$ of $S$-schemes and an element $g\in G(T')$ such that $f^{*}u = g\cdot u^{'}$.
	\end{definition}
	
	\begin{rmk}
		Let $G$ be a group scheme. The stack $J^nG$ is representable by a group scheme, which is isomorphic to $W_{n*}G$. Indeed, by Lemma \ref{Sheaf}, the functor $J^n$ coincides with the functor $W_{n*}$. Since $W_{n*}$ is defined via the functor of points, $W_{n*}G$ is a sheaf of abelian groups. Finally, by \cite[Theorem~ 12.1]{article3}, $W_{n*}G$ is representable by a scheme.
	\end{rmk}

	\begin{Lem}
		
		Let $G$ be group scheme and $U$ be an $S$-scheme. If $U$ has a $G$-action, then there is an induced action of $W_{n*}G$ on $W_{n*}U$ of $S$-schemes.
	\end{Lem}
	\begin{proof}
		By Yoneda lemma it is enough to define the action via functor of points.
		Let $T$ be an $S$-scheme, Consider $W_{n*}G(T)$ and  $W_{n*}U(T)$. Note that there is an action of $G(W_{n}^{*}(T))$ on $U(W_{n}^{*}(T))$. Let $g\in J^{n}G(T)$ and $u\in J^{n}U(T)$, by adjunction there is a unique $\tilde{g}\in G(W_{n}^{*}(T))$ and $\tilde{u}\in U(W_{n}^{*}(T))$. Now define $g\cdot u $ as the adjoint image of $ \tilde{g}\cdot \tilde{u}$.
	\end{proof}
	
	\begin{Prpn}\label{lem:QPRE}
		Let $G$ be a group scheme and $U$ be a scheme with a $G$ action, both defined over a base scheme $S$, then 
		\begin{center}
			$J^{n}[U/G]^{pre}\simeq [W_{n*}U/ W_{n*}G]^{pre} $
		\end{center}
		as prestacks, where the action of $W_{n*}G$ on  $W_{n*}U$ is as defined before.
	\end{Prpn}
	\begin{proof}
		We define a morphism of prestacks,
		\begin{center}
			$\Phi: J^{n}[U/G]^{pre}\rightarrow [W_{n*}U/ W_{n*}G]^{pre} $
		\end{center}
		Let $T$ be an $S$-scheme and $f:T'\rightarrow T $ be an S-morphism of schemes. An object over $T$ in the prestack $J^{n}[U/G]^{pre}$ is of the form $(W_{n}^{*}T,u)$ where $u:W_{n}^{*}T\rightarrow U$. By adjunction there is a unique $\tilde{u}: T\rightarrow W_{n*}U$ (we use the notation\hspace{1mm} $\tilde{}$\hspace{1mm} to denote the image under the adjunction throughout the proof) . Now we can define
		\begin{center}
			$\Phi((W_{n}^{*}T,u)): = (T,\tilde{u})$
		\end{center}
		A morphism $(W_{n}^{*}T',u')\overset{(f,g)}{\rightarrow}(W_{n}^{*}T,u)$ over $f:T'\rightarrow T$, in the category $J^{n}[U/G]^{pre}$ is the collection of data $W_{n}^{*}(f):W_{n}^{*}T'\rightarrow W_{n}^{*}T$ of $S$-schemes and an element $g\in G(W_{n}^{*}(T'))$ such that $g\cdot u' = u \circ W_{n}^{*}(f)$. Declare
		\begin{center}
			$\Phi((W_{n}^{*}T',u')\overset{(f,g)}{\rightarrow}(W_{n}^{*}T,u)) := (T',\tilde{u'})\overset{(f,\tilde{g})}{\rightarrow}(T,\tilde{u})$
		\end{center}
		which is obtained from $(f,g)$ via adjunction. Explicitly, by the definition of the action, $g\cdot u'$ is the adjoint image of $ \tilde{g}\cdot \tilde{u'}$, by adjunction again $ u \circ W_{n}^{*}(f)$ is mapped to $\tilde{u}\circ f$, which gives $\tilde{g}\cdot \tilde{u'}=\tilde{u}\circ f$. By using the same technique via adjunction we have a functor back, whose composition with $\Phi$ gives the identity functor in the respective categories.

	\end{proof}
	
	\begin{definition}[Quotient Stack]\cite[Definition.~2.4.15]{textbook2}
		The quotient stack $[U/G]$ is the category over $Sch/S$ consisting of diagrams\\
		\begin{center}
			$\xymatrix{ P \ar[d] \ar[r] &U\\ T }$
		\end{center}
		where $P\rightarrow T$ is a principal $G$-bundle and $P\rightarrow U$ is a $G$-equivariant morphism of $S$-schemes. A morphism $\left(T^{\prime} \leftarrow P^{\prime} \rightarrow U\right) \rightarrow(T \leftarrow P \rightarrow U)$ consists of a morphism $T^{\prime} \rightarrow T$ and a $G$-equivariant morphism $P^{\prime} \rightarrow P$ of schemes such that the diagram\\   
		\begin{center}
			$\xymatrix{P' \ar[r] \ar[d] \ar@/^1pc/[rr] & P \ar[r] \ar[d] & U \\
				T' \ar[r] & T \ar[u]}$	
		\end{center}	
		is commutative and the left square is cartesian.
	\end{definition}

	\begin{rmk}
		From \cite[8.1.12]{textbook1} we know that if $X$ is an algebraic space and $G/S$ is a smooth group scheme which acts on $X$, then the quotient stack $[X/G]$ is an algebraic stack.
	\end{rmk}
	
	\begin{rmk}\label{lem:QSTA}
		From \cite[2.5.21(b)]{textbook2}, if $G$ is a smooth affine group scheme acting on a scheme $U$ both defined over a scheme $S$, then the stack $[U/G]$ is isomorphic to the stackification of $[U/G]^{pre}$ over $(Sch/S)_{\acute{e}t}$.
	\end{rmk}
	
	\begin{Thm}\label{lem:jetquot}
		Let $G\rightarrow S$ be a smooth affine group scheme acting on a scheme $U$ over $S$, then $(J^{n}[U/G]^{pre})^{st}$ and $ [W_{n*}U/ W_{n*}G]$ are ismorphic as stacks.                                                               
	\end{Thm}
	
	\begin{proof}
		By Proposition \ref{lem:QPRE}, we have 	$J^{n}[U/G]^{pre}\simeq [W_{n*}U/ W_{n*}G]^{pre} $. But from Remark \ref{lem:QSTA} the stackification of $[U/G]^{pre}$ over $(Sch/S)_{\acute{e}t}$ is isomorphic to the quotient stack $[U/G]$, Now considering $(J^{n}[U/G]^{pre})^{st}$, which isomorphic to $([W_{n*}U/W_{n*}G]^{pre})^{st}$, which in turn is isomorphic to $[W_{n*}U/W_{n*}G]$ since $W_{n*}G$ is also a smooth affine group scheme over $S$.
	\end{proof}
	\begin{definition}[Classifying stacks]\cite[Definition.~2.4.14]{textbook2}
		The classifying stack $BG$ of a smooth affine group scheme $G\rightarrow S$ is the category over $Sch/S$ whose objects are principal $G$-bundles $P\rightarrow T$ and a morphism $(P'\rightarrow T')\rightarrow (P\rightarrow T)$ is the data of a $G$-equivariant morphism $P'\rightarrow P$ such that 
		
		\[
		\begin{tikzcd}
			P' \arrow[r] \arrow[d] & P \arrow[d] \\
			T' \arrow[r]            & T
		\end{tikzcd}
		\]
		is cartesian.
		
	\end{definition}	
	
	\begin{Prpn}
		Let $G$ be a smooth affine group scheme over $S$, then the stackification of the prestack $J^{n}[S/G]^{pre}$ denoted by  $(J^{n}[S/G]^{pre})^{st}$ is isomorphic to $[S/W_{n*}G]$  .
	\end{Prpn}
	
	\begin{proof}
		In Theorem \ref{lem:jetquot} we have shown that $(J^{n}[S/G]^{pre})^{st} \simeq [W_{n*}S/W_{n*}G]$ but $W_{n*}S=S$, hence we have
		\begin{center}
			$(J^{n}[S/G]^{pre})^{st} \simeq  [S/W_{n*}G]$
		\end{center} 
		as stacks.
	\end{proof}
	Now, we can consolidate the above in the form below.
	\begin{Prpn}Let $T$ be any $S$-scheme and let any object of $[U/G](W_{n}^{*}T)$ be étale locally trivial on open cover of the form $\{W_{n}^{*}V_{i}\}_{i\in I}$ for $W_{n}^{*}T$ , then the stack $(J^{n}[U/G]^{pre})^{st}$ is isomorphic to the stack $J^{n}[U/G]$.
	\end{Prpn}
	
	\begin{proof}
		Given any $S$-scheme $T$, consider the groupoid $J^{n}[U/G](T)$, which is same as the groupoid of principal $G$-bundles over $W_{n}^{*}T$, which is given by $[U/G](W_{n}^{*}T)$. Let $P\in Ob([U/G](W_{n}^{*}T))$, since it is étale locally trivial on open sets of the form $\{W_{n}^{*}V_{i}\}_{i\in I}$, the descent data $P_{i}$ of $P$ is in $[U/G]^{pre}(W_{n}^{*}V_{i})$, which is same as $[W_{n*}U/W_{n*}G]^{pre}$, which means $P$ is in the stackification $[W_{n*}U/W_{n*}G]$, hence the stack $(J^{n}[U/G]^{pre})^{st}$ is isomorphic to the stack $J^{n}[U/G]$.
	\end{proof}
	
	\section{Properties of stacks and morphisms}
	The main objective of this section is to prove that $J^{n}$ respects algebraic and Deligne-Mumford stacks. We also establish an adjunction of the Witt vector functor and $J^{n}$ functor and we thereby give a proof of the representability of jet space functor applied to a representable morphism of stack.
	\begin{Lem}\label{lem:adjunction}
		Let $T$ be a scheme, $\mathscr{X}$ be a prestack and $f:T\rightarrow J^{n}\mathscr{X}$ be a morphism of prestacks, then there exist a morphism $\tilde{f}: W_{n}^{*}T\rightarrow \mathscr{X}$ such that $ {J^{n}(\tilde{f})} \circ \eta \Rightarrow f$. (Here $\Rightarrow $ denote a 2-isomorphism of a morphism of stacks.)
		\[
		\xymatrix@C=1cm@!C{
			T \ar[r]^{\eta} & J^n W_n^* T \ar[r]^{J^n(\tilde{f})} & J^n \mathscr{X}
		}
		\]

		where $\eta $ is the unit map of the adjuction between $J^{n}$ and $W_{n}^{*}$.
	\end{Lem}
	
	\begin{proof}
		The primary objective is to construct a functor $\tilde{f}$ which is a candidate for the above required properties. Note that all the objects considered in this proof are stacks ie. If $T$ is a scheme then we are working with a stack $\mathscr{X}_{T}$ representing the scheme $T$. \\
		Consider $Id :T\rightarrow T $ be the identity morphism of schemes then $(T,Id)$ is an object in the stack associated to the scheme $T$, ie.
		\begin{center}
			$\xymatrix{(T,Id) \ar[d]^{p} \\ T}$
		\end{center}
		$(T,Id)$ is an object over $T$, by the definition of the morphism between stacks, we have $f((T,Id))\in J^{n}\mathscr{X}(T)$, but by definition $f((T,Id))\in \mathscr{X}(W_{n}^{*}T)$. Now consider an object $(U,a)\in Ob(W_{n}^{*}T)$, in the stack $\mathscr{X}$ we have the following commutative square
		\begin{center}
			$	\xymatrix{f((T,Id))|_{U}\ar@{|.>}[r]\ar@{|->}[d] &f((T,Id))\ar@{|->}[d]\\ U \ar[r]^a &W_{n}^{*}T} $
		\end{center}
		since the pullbacks are unique upto unique isomorphism, we define
		\begin{center}
			$\tilde{f}((U,a)) = f((T,Id))|_{U}$
		\end{center}
		let $(V,b)\rightarrow (U,a)$ be a morphism in $W_{n}^{*}T$, consisting of the following data $g:V\rightarrow U$, a morphism of schemes and $b= g^{*}a$. Now we have to determine the image of the above morphism. We have the following data in the category $\mathscr{X}$
		\begin{center}
			$	\xymatrix{ f((T,Id))|_{V}\ar@{|.>}[r]\ar@{|->}[d]^{p}\ar@/^1pc/[rr] &f((T,Id))|_{U}\ar@{->}[r]\ar@{->}[d]^{p} &f((T,Id))\ar@{|->}[d]^{p}\\ V \ar@{->}[r]^{g}&U \ar[r]^a &W_{n}^{*}T} $
		\end{center}   
		By axiom(2) of prestacks the dotted arrow in the above diagram is unique. we denote the morphism as $\tilde{f}((g,g^{*}))$. That is,
		\begin{center}
			$\tilde{f}((V,u)\overset{(g,g^{*})}\rightarrow (U,a)) :=\tilde{f}((g,g^{*})) $
		\end{center}
		The morphism $\eta : T\rightarrow J^{n}W_{n}^{*}T$, which is the unit map of the adjunction of the functors $J^{n}$ and $W_{n}^{*}$, which is represented as a morphism of stacks as; Let $(U,h)\in Ob(T)$.
		\begin{center}
			$\eta ((U,h)) = (W_{n}^{*}U,W_{n}^{*}(h))$
		\end{center}
		Now we will unwind the above definition, $(U,h)\in Ob(T)$ is the data consisting of the morphism of schemes $h:U\rightarrow T$. For defining a functor $\eta $ we have to produce an object in $J^{n}W_{n}^{*}T(U)$, which is the same as giving an element in  $W_{n}^{*}T(W_{n}^{*}U)$. By observing that the morphism $h$ induces a morphism functorially, which is given by $W_{n}^{*}U \overset{W_{n}^{*}(h)}\rightarrow W_{n}^{*}T$, which is same as producing an element $(W_{n}^{*}U,W_{n}^{*}(h))\in Ob(J^{n}W_{n}^{*}T(U))$. Now suppose $(V,t)\rightarrow (U,h)$ be a morphism in $T$ with the data given by $V\overset{s}\rightarrow U $ such that $t = s^{*}h$. By the above definition of mapping of objects $(V,t)\mapsto (W_{n}^{*}V,W_{n}^{*}(t))$ similarly $(U,h)\mapsto (W_{n}^{*}U,W_{n}^{*}(h))$. Now we can evidently define the map to be,
		\begin{center}
			$\eta((V,t)\overset{s}\rightarrow (U,h)) =  (W_{n}^{*}V,W_{n}^{*}(t))\overset{W_{n}^{*}(s)}\rightarrow (W_{n}^{*}U,W_{n}^{*}(h))$.
		\end{center}
		
		We have a functor $J^{n}(\tilde{f}): J^{n}W_{n}^{*}T\rightarrow J^{n}\mathscr{X}$, induced from the morphism $\tilde{f}:W_{n}^{*}T\rightarrow \mathscr{X}$, which is as follows\\
		Let $U$ be a scheme, by the definition of the functor $J^{n}$, an object over $U$ is given by $(W_{n}^{*}U,l)$ where $l:W_{n}^{*}U\rightarrow W_{n}^{*}T$, so  we have  $(W_{n}^{*}U,l)\in Ob(J^{n}W_{n}^{*}T(U)) =Ob(W_{n}^{*}T(W_{n}^{*}U)) $, but $\tilde{f}:W_{n}^{*}T\rightarrow \mathscr{X}$, which maps objects in $W_{n}^{*}T(W_{n}^{*}U)$ to $\mathscr{X}(W_{n}^{*}U) = J^{n}\mathscr{X}(U)$. Hence we can define,
		\begin{center}
			$J^{n}(\tilde{f})((W_{n}^{*}U,l)) = \tilde{f}((W_{n}^{*}U,l))$
		\end{center}
		Now let $(W_{n}^{*}V,k)\overset{( W_{n}^{*}(s), W_{n}^{*}(s)^{*})}\longrightarrow (W_{n}^{*}U,l)$ be a morphsim in $J^{n}W_{n}^{*}T$, which is given by the following data 
		
		\begin{center}
			$	\xymatrix{(W_{n}^{*}V,k) \ar@{|->}[d]^{p} \ar@{->}[r]  & (W_{n}^{*}U,l) \ar@{|->}[d]^{p}\\W_{n}^{*}V\ar@{->}[r]^{W_{n}^{*}(s)}\ar@{|->}[d]^{p'} &W_{n}^{*}U\ar@{|->}[d]^{p'}\\ V \ar[r]^s &U} $
		\end{center}
		where  $k = W_{n}^{*}(s)^{*}l$. The morphism $\tilde{f}((W_{n}^{*}V,k)\overset{( W_{n}^{*}(s), W_{n}^{*}(s)^{*})}\longrightarrow (W_{n}^{*}U,l))$ is the following data in the stack $J^{n}\mathscr{X}$
		
		\begin{center}
			$	\xymatrix{\tilde{f}((W_{n}^{*}V,k)) \ar@{|->}[d]^{p} \ar@{->}[r]  & \tilde{f}((W_{n}^{*}U,l)) \ar@{|->}[d]^{p}\\W_{n}^{*}V\ar@{->}[r]^{W_{n}^{*}(s)}\ar@{|->}[d]^{p'} &W_{n}^{*}U\ar@{|->}[d]^{p'}\\ V \ar[r]^s &U} $
		\end{center}
		
		Now we will look at the composite functor $\eta \circ J^{n}(\tilde{f})$\\
		Let $(U,h)\in Ob(T)$, via composition we have below sequence of mapping
		
		\begin{align}
			(U,h) \mapsto (W_{n}^{*}U,W_{n}^{*}(h)) \mapsto \tilde{f}((W_{n}^{*}U,W_{n}^{*}(h)))= (W_{n}^{*}(h))^{*}f(T,Id)    
		\end{align}
		
		Note that $f((U,h)) \in J^{n}\mathscr{X}(U)$, also there is a morphism in $T$ given by $(U,h)\rightarrow(T,Id)$, which gives the following diagram
		
		\begin{center}
			$	\xymatrix{f((U,h)) \ar@{|->}[d]^{p} \ar@{->}[r]  & f((T,id)) \ar@{|->}[d]^{p}\\W_{n}^{*}U\ar@{->}[r]^{W_{n}^{*}(h)}\ar@{|->}[d]^{p'} &W_{n}^{*}T\ar@{|->}[d]^{p'}\\ U \ar[r]^h &T} $
		\end{center}
		
		Consider the below diagram 
		
		\begin{center}
			\[
			\xymatrix@C=3em@R=3em{
				(W_n^*(h))^* f((T,\mathrm{Id})) \ar@{-->}[r] \ar@{|->}[d]^{p} \ar@/^1.5pc/[rr] & f((U,h)) \ar[r] \ar[d]^{p} & f((T,\mathrm{Id})) \ar@{|->}[d]^{p} \\
				W_n^* U \ar[r]^{\mathrm{Id}} \ar[d] & W_n^* U \ar[r]^{W_n^*(h)} \ar[d] & W_n^* T \ar[d] \\
				U \ar[r]^{\mathrm{Id}} & U \ar[r]^h & T
			}
			\]
		\end{center} 
		
		The dotted arrow in the above diagram is a unique isomorphism, giving the fact that $ \eta \circ J^{n}(\tilde{f})$ and $f$ are 2-isomorphic.
			
	\end{proof}
	
	 Now we will give the proof of Theorem \ref{intro_ad} as a consequence of Lemma \ref{lem:adjunction}.
	\begin{Cor}
		Let $T$ be a scheme and $\mathscr{X}$ be a prestack over $(Sch/S)_{\acute{e}t}$, then
		
		\begin{center}
			$Hom(T,J^{n}\mathscr{X})\simeq Hom(W_{n}^{*}T,\mathscr{X})$
		\end{center}
		upto natural equivalence of morphisms of prestacks.
	\end{Cor}
	\begin{proof}
		Consider the map $\Psi :Hom(T,J^{n}\mathscr{X}) \rightarrow Hom(W_{n}^{*}T,\mathscr{X}) $, defined as
		\begin{center}
			$\Psi (f): = \tilde{f}$ , $\forall f\in Hom(T,J^{n}\mathscr{X})$
		\end{center}
		where $\tilde{f}$ is the same functor constructed in Lemma \ref{lem:adjunction}
		\vspace{0.3cm} 
		\\Now consider the map $\Theta : Hom(W_{n}^{*}T,\mathscr{X}) \rightarrow Hom(T,J^{n}\mathscr{X})$, defined as
		\begin{center}
			$\Theta (g): = \eta \circ J^{n}(g)$, $\forall g\in Hom(W_{n}^{*}T,\mathscr{X})$. 
		\end{center}
		Now it remains to show that $\widetilde{\eta \circ J^{n}(g)}$, defined in Lemma \ref{lem:adjunction}  is equivalent to $g$ as  morphisms of stacks. Consider the below diagram
		
		\begin{center}
			\[
			\xymatrix@C=3em@R=3em{
				g(U,a) \ar@{.>}[r] \ar@/^1pc/[rr]\ar[d] & a^{*}g(W_{n}^{*}T,Id) \ar[r] \ar[d] & g(W_{n}^{*}T,Id) \ar[d] \\
				U \ar[r]^{\mathrm{Id}} & U \ar[r]^a & W_{n}^{*}T
			}
			\]
		\end{center} 
		here the arrow $g(U,a) \rightarrow g(W_{n}^{*}T,Id)$ is obtained from the arrow $(U,a)\rightarrow (W_{n}^{*}T,Id)$ functorially, by using $g$ and looking $W_{n}^{*}T$ as a stack. Now $\widetilde{\eta \circ J^{n}(g)}(T,id)$ is by definition of $\eta$, we have
		\begin{center}
			$\widetilde{\eta \circ J^{n}(g)}(T,id) =  g(W_{n}^{*}T,Id)$
		\end{center} 
		Now from the diagram above, it is evident that $g(U,a)\rightarrow a^{*}g(W_{n}^{*}T,Id)$ is an isomorphism, giving the required adjunction up to 2-isomorphism of morphisms of stacks.	
	\end{proof}
	\begin{rmk}
		In the case of sheaves, the pullbacks are unique and we have bijection as sets, which is same as the adjunction of functors $W_{n}^{*}$ and $W_{n*}$.
	\end{rmk}
	
	\begin{Prpn}\label{rep_mor}
		Let $f:\mathscr{X} \rightarrow \mathscr{Y}$ be a representable morphism of prestacks, then $J^{n}(f):J^{n}\mathscr{X} \rightarrow J^{n}\mathscr{Y}$ is representable.
	\end{Prpn}
	
	\begin{proof}
		For the representability of $J^{n}(f)$, we have to show that for any given morphism $T \xrightarrow{h} J^{n}\mathscr{Y}$ from a scheme, the fibre product $J^{n}\mathscr{X}\times _{J^{n}\mathscr{Y}} T $ is an algebraic space. By Lemma \ref{lem:adjunction}, for the morphism $h: T \rightarrow J^{n}\mathscr{Y}$, there exist a morphism $\tilde{h} :W_{n}^{*}T \rightarrow \mathscr{Y}$ such that 
		\[
		\xymatrix@C=1cm@!C{
			T \ar[r]^{\eta} & J^n W_n^* T \ar[r]^{J^n(\tilde{h})} & J^n \mathscr{Y}
		}
		\]
		the composition of above morphism gives, $J^{n}(\tilde{h})\circ \eta \Rightarrow h$. Now since the morphism $f$ is representable, we have
		\begin{center}
			$\mathscr{X}\times _{f,\mathscr{Y},\tilde{h}}W_{n}^{*}T$
		\end{center}
		representable by an algebraic space. Now applying the $J^{n}$ functor, we have
		\begin{center}
			$J^{n}(\mathscr{X}\times _{f,\mathscr{Y},\tilde{h}}W_{n}^{*}T) \simeq J^{n} \mathscr{X}\times _{J^{n}(f),J^{n} \mathscr{Y},J^{n} (\tilde{h})}J^{n}W_{n}^{*}T 
			$
		\end{center}
		Now taking the base change along the map $\eta : T\rightarrow J^{n}W_{n}^{*}T $
		\begin{center}
			$ J^{n} \mathscr{X}\times _{J^{n}(f),J^{n} \mathscr{Y},J^{n} (\tilde{h})}J^{n}W_{n}^{*}T\times_{id,J^{n}W_{n}^{*}T,\eta} T \simeq J^{n} \mathscr{X}\times _{J^{n}(f),J^{n} \mathscr{Y},h}T$
		\end{center}
		here the last isomorphism is due to the fact that $J^{n}(\tilde{h})\circ \eta \Rightarrow h $ and applying properties of base change.	 
	\end{proof}
	
	\begin{Cor}\label{rep-diagonal}
		Let $\mathscr{X}$ be an algebraic stack, then the diagonal $\Delta _{J^{n}\mathscr{X}}$ is representable.
	\end{Cor}
	
	\begin{proof}
		Note that the diagonal of $J^{n}\mathscr{X}$ is obtained by applying $J^{n}$ to the diagonal of $\mathscr{X}$, since the diagonal of the algebraic stack $\mathscr{X}$ is representable, using Proposition \ref{rep_mor}, we have the representability of the diagonal $\Delta _{J^{n}\mathscr{X}}$ of $J^{n}\mathscr{X}$. 
	\end{proof}

	\begin{Prpn}\label{lem:surjection}
		Let $f:U\rightarrow \mathscr{X}$ be a smooth (resp. étale) surjective morphism from a scheme to  a stack, then $J^{n}(f):J^{n}U\rightarrow J^{n}\mathscr{X}$ is smooth (resp. étale) surjective morphism.
	\end{Prpn}
	
	\begin{proof}
		Let $T$ be a scheme and $g: T\rightarrow J^{n}\mathscr{X}$ be a morphism, our aim is to show that the projection morphism $J^{n}U\times _{J^{n}\mathscr{X}}T\xrightarrow{p_{2}} T$ is smooth(resp. étale) surjective morphism.\\
		Note that by Corollary~\ref{rep-diagonal} the diagonal morphism $\Delta _{J^{n}\mathscr{X}}$ is representable, which gives that $J^{n}U\times _{J^{n}\mathscr{X}}T$ is an algebraic space.
		Given $g:T\rightarrow J^{n}\mathscr{X}$ by Lemma~\ref{lem:adjunction}, there exists $\tilde{g}: W_{n}^{*}T\rightarrow \mathscr{X}$ such that $J^{n}(\tilde{g})\circ \eta \Rightarrow g$. Now, the map
		\begin{center}
			$U\times _{f,\mathscr{X},\tilde{g}}W_{n}^{*}T\xrightarrow{p_{2}} W_{n}^{*}T$
		\end{center}
		is smooth (resp.étale) surjective, using the fact that $f$ has the prior property.
		Applying the functor $J^{n}$ to the above map, we have
		\begin{center}
			$J^{n}(U\times _{f,\mathscr{X},\tilde{g}}W_{n}^{*}T)\xrightarrow{J^{n}(p_{2})} J^{n}W_{n}^{*}T$
		\end{center}
		Using the fact that $J^{n}$ behave well with the fiber products, we have
		
		\begin{center}
			$J^{n}U\times _{J^{n}(f),J^{n}\mathscr{X},J^{n}(\tilde{g})}J^{n}W_{n}^{*}T\xrightarrow{p_{2}} J^{n}W_{n}^{*}T$
		\end{center}
		From \cite[Proposition  13.3(h)]{article3} the functor $J^{n}$ preserves smooth and surjective map but étale morphisms are also smooth hence the last arrow is also smooth (resp.étale) surjective. Now making the base change via the co-unit map, we have
		\begin{center}
			$J^{n}U\times _{J^{n}(f),J^{n}\mathscr{X},J^{n}(\tilde{g})}(J^{n}W_{n}^{*}T\times _{Id,J^{n}W_{n}^{*}T,\eta}T)\xrightarrow{p_{2}} J^{n}W_{n}^{*}T\times _{Id,J^{n}W_{n}^{*}T,\eta}T$
		\end{center}   
		which is same as
		\begin{center}
			$J^{n}U\times _{J^{n}(f),J^{n}\mathscr{X},g} T\xrightarrow{p_{2}} T$
		\end{center}
		since $J^{n}(\tilde{g})\circ \eta $ and $ g$ are 2-isomorphic, we have $p_{2}$ as smooth surjective (resp. étale surjective).
		
	\end{proof}
	
	\begin{Thm}
		Let $\mathscr{X}$ be an algebraic stack (resp. Deligne Mumford stack), then $J^{n}\mathscr{X}$ is an algebraic stack (resp. Deligne Mumford stack).
	\end{Thm}
	
	\begin{proof}
		By Corollary \ref{rep-diagonal}, the diagonal $\Delta_{J^{n}\mathscr{X}}: J^{n}\mathscr{X}\times _{S}J^{n}\mathscr{X}$ is representable and by Proposition \ref{lem:surjection}, there is a smooth (resp. étale surjection) from a scheme, which gives that $J^{n}\mathscr{X}$ is algebraic (resp. Deligne Mumford) stack.

	\end{proof}
	
	\begin{definition}[quasi-compact]\cite[Lemma 100.62 (3)]{Stk_pro_2}
		An algebraic stack $\mathscr{X}$ is quasi-compact if there is a surjective smooth morphism $U\rightarrow \mathscr{X}$ with $U$ a quasi compact scheme.
	\end{definition}
	
	\begin{definition}[quasi-separated]\cite[Definition 8.2.12]{textbook1}
		Let $f:\mathscr{X}\rightarrow \mathscr{Y}$ be a morphism of algebraic stacks over $S$ , and let
		\begin{center}
			$\Delta_{\mathscr{X}/\mathscr{Y}}:\mathscr{X}\rightarrow \mathscr{X} \times_{\mathscr{Y}}\mathscr{X}$
		\end{center}
		be the diagonal morphism.
		$f$ is quasi-separated if the diagonal $\Delta_{\mathscr{X}/\mathscr{Y}} $ is quasi-compact and quasi-separated. In the case where $\mathscr{Y}=S$ and $f$ is the structure morphism, we say that $\mathscr{X}$ is quasi-separated.
		
	\end{definition}

	\begin{definition}[locally Noetherian stack] \cite[Definition 4.3.7]{textbook2}  
		Let $\mathscr{X}$ be an algebraic stack. We say that  $\mathscr{X}$ is locally Noetherian if there exists a smooth presentation $U\rightarrow \mathscr{X}$ such that $U$ is locally Noetherian.
		
	\end{definition}
	
	\begin{definition}[Noetherian Stack]\cite[Definition 0510]{Stk_pro_1}
		Let $\mathscr{X}$ be an algebraic stack. We say $\mathscr{X}$ is Noetherian if it is quasi-compact, quasi-separated and locally Noetherian.
	\end{definition}

	\begin{Prpn}\label{qcqs-stacks}
		Let $\mathscr{X}/S$ be an algebraic stack, then we have the following
		\begin{enumerate}[label=(\roman*)]
			\item If $\mathscr{X}$ is quasi-compact, then $J^{n}\mathscr{X}$ is quasi-compact.
			\item Let $f:\mathscr{X}\rightarrow \mathscr{Y}$ be a quasi-compact morphism of algebraic stack, which is representable then $J^{n}(f):J^{n}\mathscr{X}\rightarrow J^{n}\mathscr{Y}$ is quasi-compact.
			\item \label{qsprtd} If $\mathscr{X}$ is quasi-separated, then $J^{n}\mathscr{X}$ is quasi-separated.

			\item Let $f: \mathscr{X}\rightarrow \mathscr{Y}$ be a quasi-separated morphism, then $J^n(f): J^n\mathscr{X}\rightarrow J^n\mathscr{Y}$ is quasi-separated.
			\item Let $S$ be locally Noetherian scheme, and $\mathscr{X}$ be locally of finite type over $S$, then $J^{n}\mathscr{X}$ is locally Noetherian.
		\end{enumerate} 
	\end{Prpn}
	\begin{proof}
		\begin{enumerate}[label=(\roman*)]
			\item As $\mathscr{X}$ is quasi-compact, there exists a surjective smooth morphism $U\rightarrow \mathscr{X}$ where $U$ is a quasi-compact scheme. From \cite[Proposition 11.10 (a)]{article3} $J^{n}U$ is quasi-compact and using Proposition \ref{lem:surjection}, $J^{n}U\rightarrow J^{n}\mathscr{X}$ is smooth and surjective. Hence $J^{n}\mathscr{X}$ is quasi-compact.
			\item Let $g: V\rightarrow J^{n}\mathscr{Y}$ be a morphism such that $V$ is a scheme, we will show that $J^{n}\mathscr{X}\times _{J^{n}\mathscr{Y}}V\rightarrow V$ is quasi-compact. By Lemma \ref{lem:adjunction} there is a map $\tilde{g}:W_{n}^{*}V\rightarrow \mathscr{Y}$ such that $J^{n}(\tilde{g})\circ \eta \Rightarrow g$. Since $f$ is representable and quasi-compact, the map $\mathscr{X}\times _{\mathscr{Y}}W_{n}^{*}V \rightarrow W_{n}^{*}V$ is a quasi-compact morphism of algebraic spaces. By \cite[Proposition 11.10 (c)]{article7} $J^{n}$ preserves quasi-compact morphism at the level of algebraic spaces, which gives the quasi-compact map $J^{n}\mathscr{X}\times _{J^{n}\mathscr{Y}}J^{n}W_{n}^{*}V \rightarrow J^{n}W_{n}^{*}V$, now by taking the base change along the map $\eta :V \rightarrow J^{n}W_{n}^{*}V$, we have the quasi-compact map $J^{n}\mathscr{X}\times _{J^{n}\mathscr{Y}}V\rightarrow V$ using the fact that $J^{n}(\tilde{g})\circ \eta \Rightarrow g$.
			\label{quasi-compact}
			\item For $J^{n}\mathscr{X}$ to be quasi-separated, we have to show that its diagonal and second diagonal (ie, $\mathscr{X}\rightarrow \mathscr{X}\times _{\mathscr{X}\times _{S}\mathscr{X}}\mathscr{X}$) is quasi-compact. The diagonal of $J^{n}\mathscr{X}$ is $J^{n}$ applied to the diagonal, same is true for the second diagonal as well. Now apply Proposition \ref{qcqs-stacks}(\ref{quasi-compact}) for diagonal and second diagonal. 
			\item proof is same as \ref{qsprtd} with the diagonal replaced with the diagonal with respect to $f$.
			\item Since $\mathscr{X}$ is locally of finite type, there is a smooth surjection $U\rightarrow \mathscr{X}$ such that $U$ is a scheme locally of finite type. From \cite[Proposition 13.3 (c)]{article3} $J^{n}U$ is locally of finite type. Now using the fact that a scheme locally of finite type over a Noetherian base is locallly Noetherian \cite[Lemma 29.15.6]{Stk_pro_3}, $J^{n}U$ is locally Noetherian. Hence we have a smooth surjection $ J^{n}U\rightarrow J^{n}\mathscr{X}$ where $J^{n}U$ is locally Noetherian over $S$.
		\end{enumerate}
	\end{proof}
	
	\begin{definition}[$\underline{Isom}_{\mathscr{X}}$ sheaf] \cite[3.4.7]{textbook1}
		Let $\mathscr{X}$ be a prestack, $U \in (Sch/S)$ and $x,x'\in \mathscr{X}(U)$, We define a presheaf
		\begin{center}
			$\underline{Isom}_{\mathscr{X}}(x,x'):(Sch/U)^{op}\rightarrow Set $
		\end{center}
		as \\
		For any morphism $V\xrightarrow{f} U$, consider pullbacks $f^{*}x,f^{*}x'$, we set
		\begin{center}
			$\underline{Isom}_{\mathscr{X}}(x,x')(f:V\rightarrow U): = Isom_{\mathscr{X}(V)}(f^{*}x,f^{*}x')$
		\end{center}
		For a composition, 
		\begin{center}
			$Z\xrightarrow{g} V \xrightarrow{f}U$
		\end{center}
		the pullback $(fg)^{*}x$ is a pullback along  $g$ of $f^{*}x$ (\text{resp.} $f^{*}x')$, and therefore there is a canonical map
		\begin{center}
			$g^{*}: Isom_{\mathscr{X}(V)}(V\xrightarrow{f} U) \rightarrow Isom_{\mathscr{X}(W)}(W\xrightarrow{f\circ g} U) $
		\end{center}
		
	\end{definition}
	
	\begin{rmk}
		Let $\mathscr{X}$ be a stack and consider $J^{n}\mathscr{X}$, for any $S$-scheme $U$ and any $x,x'\in J^{n}\mathscr{X}(U)$, the sheaf $\underline{Isom}_{J^{n}\mathscr{X}}(x,x')$ on $(Sch/U)$ is a sheaf. Now we consider, 
		\begin{align}
			\underline{\mathrm{Isom}}_{J^{n}\mathscr{X}}(x, x')(f\colon V \rightarrow U) 
			:&= \mathrm{Isom}_{J^{n}\mathscr{X}(V)}(f^{*}x, f^{*}x') \\
			&= \mathrm{Isom}_{\mathscr{X}(W_{n}^{*}V)}((W_{n}^{*}f)^{*}x, (W_{n}^{*}f)^{*}x') 
			\label{Isomsheaf1}
		\end{align}
		where the second equality holds by definition of the jet functor.\\
		Note that $x,x'\in \mathscr{X}(W_{n}^{*}U)$, we now look at the structure of $\underline{\mathrm{Isom}}_{J^{n}\mathscr{X}}(x, x') $ compared to $\underline{\mathrm{Isom}}_{\mathscr{X}}(x, x')$. Now evaluating $\underline{\mathrm{Isom}}_{\mathscr{X}}(x, x') $ on $W_{n}^{*}V$ for any $V$ over $U$. We call the resulting presheaf as $J^{n}\underline{\mathrm{Isom}}_{\mathscr{X}}(x, x')$.
		\begin{align}
			J^{n}\underline{\mathrm{Isom}}_{\mathscr{X}}(x, x')(f\colon V \rightarrow U) 
			:&=\underline{\mathrm{Isom}}_{\mathscr{X}}(x,x')(W_{n}^{*}f:W_{n}^{*}
			V\rightarrow W_{n}^{*}U) \\
			&= \mathrm{Isom}_{\mathscr{X}(W_{n}^{*}(V))}((W_{n}^{*}f)^{*}x, (W_{n}^{*}f)^{*}x') 
			\label{Isomsheaf2}
		\end{align}
		From \ref{Isomsheaf1} and \ref{Isomsheaf2}, we see that the sheaves $J^{n}\underline{\mathrm{Isom}}_{\mathscr{X}}(x, x')$ and $\underline{\mathrm{Isom}}_{J^{n}\mathscr{X}}(x, x')$ are isomorphic.  
	\end{rmk}
	\begin{Prpn}
		Let $\mathscr{X}$ be a stack and the diagonal $\Delta _{\mathscr{X}}$ be representable, then the sheaf $\underline{\mathrm{Isom}}_{J^{n}\mathscr{X}}(x, x')$ defined above is an algebraic space for any $U\in (Sch/S)$ and $x,x'\in J^{n}\mathscr{X}(U)$.
	\end{Prpn}
	
	\begin{proof}
		The diagonal of $J^{n}\mathscr{X}$ is representable by Proposition \ref{rep_mor}. From \cite[Lemma 8.1.8]{textbook1} the diagonal $\Delta _{\mathscr{X}}$ of a stack $\mathscr{X}$ is reprentable if and only if for every $S$-scheme $U$ and objects $x,x' \in \mathscr{X}(U) $ the sheaf $\underline{Isom}(x,x')$ sheaf on $(Sch/U)$ is an algebraic space, which shows $\underline{\mathrm{Isom}}_{J^{n}\mathscr{X}}(x, x')$ is an algebraic space.
	\end{proof}
	
	\begin{definition}[Inertia stack]\cite[Definition 3.2.9]{textbook2} The \textit{inertia stack} of an algebraic stack $\mathscr{X}$ is the fiber product 
		\begin{equation}\label{inertia_square}
			\begin{tikzcd}
				I_{\mathscr{X}} \arrow[r] \arrow[d] 
				& \mathscr{X} \arrow[d] \\
				\mathscr{X} \arrow[r] 
				& \mathscr{X} \times \mathscr{X}
				\arrow[ul, phantom, "\square"]
			\end{tikzcd}
		\end{equation}
		
	\end{definition}
	
	\begin{Prpn}
		Let $\mathscr{X}$ be an algebraic stack, then the \textit{inertia stack} of $J^{n}\mathscr{X}$ is $J^{n}I_\mathscr{X}$.
	\end{Prpn}
	
	\begin{proof}
		Applying $J^{n}$ to the cartesian square \ref{inertia_square}, it is still cartesian and we have,
		\begin{equation}
			\begin{tikzcd}
				J^{n}I_{\mathscr{X}} \arrow[r] \arrow[d] 
				& J^{n}\mathscr{X} \arrow[d] \\
				J^{n}\mathscr{X} \arrow[r] 
				& J^{n} \mathscr{X} \times J^{n}\mathscr{X}
				\arrow[ul, phantom, "\square"]
			\end{tikzcd}
		\end{equation}
		from the above it follows that $J^{n}I_\mathscr{X} \simeq I_{J^{n}\mathscr{X}}$.
	\end{proof}
	
	\section{Morphisms between Arithmetic Jet Stacks}
	\subsection{Construction of morphisms}
	Given any scheme $X$ and $i\in \mathbb{N}$ and a collection of morphisms 
	
	\begin{center}
		$W_{n}^{*}X\xrightarrow{\phi_{X,n,i}}W_{n+i}^{*}X$
	\end{center}
	satisfying the following commutative diagram for any schemes $X,Y$
	
	\begin{align}\label{dg1}
		\xymatrix{
			W_{n}^{*}X \ar[d]_{W_{n}^{*}(f)} \ar[r]^{\phi_{X,n,i}} &W_{n+i}^{*}X\ar[d]^{W_{n+i}^{*}(f)}\\
			W_{n}^{*}Y \ar[r]^{\phi_{Y,n,i}}
			&W_{n+i}^{*}Y}     
	\end{align}
	and for any $f:X\rightarrow Y$ induce morphisms between different levels of jet stack over a given stack. We will demonstrate the general method of constructing such a morphism.\\
	
	\begin{Prpn}\label{lem:trun-frob}
		Let  $\mathscr{X}$ be a stack and $W_{n}^{*}X\xrightarrow{\phi_{X,n,i}}W_{n+i}^{*}X$ be a morphism of schemes for every scheme $X$ and satisfying the diagram \ref{dg1}, then there exist a morphism $J^{n+i}\mathscr{X}\xrightarrow{\Phi_{X,n,i}} J^{n}\mathscr{X}$, satisfying the 2-commutative diagram
        \begin{align} \label{Frob_trunc_stacks}
             \xymatrix{
          J^{n+i}\mathscr{X} \ar[r]^{\Phi_{X,n,i}} \ar[d]_{J^{n+i}(g)} & J^{n}\mathscr{X} \ar[d]^{{J^{n}(g)}} \\
         J^{n+i}\mathscr{Y} \ar[r]^{\Phi_{Y,n,i}} & J^{n}\mathscr{Y}
         }
        \end{align}
         for any $g:\mathscr{X}\rightarrow \mathscr{Y}$.
	\end{Prpn}
	
	\begin{proof}
		Let $U$ be an $S$-scheme and let $x\in J^{n+i}\mathscr{X}(U)$, consider the following diagram
		\begin{center}
			$	\xymatrix{\phi_{U,n,i}^{*}x \ar@{.>}[r]\ar@{|->}[d]&x\ar@{|->}[d]^{}\\ W_{n}^{*}U \ar[r]^{\phi_{U,n,i}} &W_{n+i}^{*}U} $
		\end{center}
		where the pull back is obtained using the fact that $\mathscr{X}$ is a prestack. Now we define the functor
		
		\begin{center}
			$\Phi_{U,n,i}(x):=\phi_{U,n,i}^{*}x$
		\end{center}
		let $y\xrightarrow{\alpha} x$ be a morphism in $J^{n+i}\mathscr{X}$, which is the data given below,
		\begin{center}
			$	\xymatrix{y \ar@{->}[r]^{\alpha}\ar@{|->}[d]&x\ar@{|->}[d]^{}\\ W_{n+i}^{*}V \ar[r]^{W_{n+i}^{*}(f)} &W_{n+i}^{*}U} $
		\end{center}
		Using the property of prestacks, as done in the previous cases, we have a unique arrow $\tilde{\alpha}$ between the pullbacks via the given map of schemes, which is illustrated by the diagram
		
		\begin{center}
			$
			\xymatrix{
				\phi_{V,n,i}^{*}y \ar@{.>}[r]^{\tilde{\alpha}} \ar@{.>}@/^1pc/[rr] \ar@{|->}[d] & \phi_{U,n,i}^{*}x \ar@{|->}[d] \ar@{.>}@/^1pc/[rr]
				& y \ar[r] \ar@{|->}[d] & x \ar@{|->}[d] \\
				W_{n}^{*}V \ar@/^2pc/[rr]_{\phi_{V,n,i}} \ar[r]^{W_{n}^{*}(f)} & W_{n}^{*}U \ar@/_1pc/[rr]_{\phi_{U,n,i}} & W_{n+i}^{*}V \ar[r]^{W_{n+i}^{*}(f)} & W_{n+i}^{*}U
			}
			$
		\end{center}
		Now we are ready to map the given morphism $\alpha$
		\begin{center}
			$\Phi_{X,n,i}(\alpha):=\tilde{\alpha}$
		\end{center}
        Now we will prove the 2-commutativity of diagram \ref{Frob_trunc_stacks}. For that we will show $\Phi_{Y,n,i}\circ J^{n+i}(g)$ and $ J^{n}(g)\circ \Phi_{X,n,i}$ are 2-isomorphic. Let $U\in (Sch/S)$ and $x\in J^{n+i}\mathscr{X}(U)$ under the map $\Phi_{Y,n,i}\circ J^{n+i}(g)$ we have the following diagram
        \begin{align*}
            \xymatrix{\phi_{U,n,i}^{*}g(x) \ar@{->}[r]\ar@{|->}[d]&g(x)\ar@{|->}[d]^{}\\ W_{n}^{*}U \ar[r]^{\phi_{U,n,i}} &W_{n+i}^{*}U}
        \end{align*}
     similarly for $J^{n}(g)\circ \Phi_{X,n,i}$ we have,
     \begin{align*}
            \xymatrix{g(\phi_{U,n,i}^{*}x) \ar@{->}[r]\ar@{|->}[d]&g(x)\ar@{|->}[d]^{}\\ W_{n}^{*}U \ar[r]^{\phi_{U,n,i}} &W_{n+i}^{*}U}
        \end{align*}
        by comparing the diagrams above $\phi_{U,n,i}^{*}g(x)$ and $g(\phi_{U,n,i}^{*}x)$ are isomorphic. We are done.
	\end{proof}
	
	\begin{rmk}
		Given any scheme $X$, $i\in \mathbb{N}$ and morphisms 
		\begin{center}
			$W_{n+i}^{*}X\xrightarrow{\tilde{\phi}_{X,n,i}}W_{n}^{*}X$
		\end{center}
		satisfying commutative diagram as in \ref{dg1}, the constructions in Proposition \ref{lem:trun-frob} follows and we have a morphism of stacks
		\begin{center}
			$ J^{n}\mathscr{X}\xrightarrow{\tilde{\Phi}_{X,n,i}} J^{n+i}\mathscr{X}$
		\end{center}
	\end{rmk}
	
	\begin{Cor}
		The trunation and  Frobenius morphisms satisfies \ref{dg1} and induces a morphism between jet arithmetic stacks as in Proposition \ref{lem:trun-frob}.
	\end{Cor}
	\begin{proof}
		Repeat of the proof of Proposition \ref{lem:trun-frob} by replacing the morphism with truncation or Frobenius.
	\end{proof}
	
	\subsection{Properties of morphisms between jet algebraic stacks.}
	
	\begin{definition}\cite[Definition~8.2.6]{textbook1}(Properties of morphisms)
		Let $f:\mathscr{X}\rightarrow \mathscr{Y}$ be a morphism of algebraic stacks over $S$. A chart for $f$ is the commutative diagram
		\begin{align}\label{morphism_stack_diagram}
			\xymatrix{
				X \ar[r]^g \ar[rd]^q \ar@/^1.5pc/[rr]^h & \mathscr{X}' \ar[r]^{f'} \ar@{->}[d]^{p'} & Y \ar@{->}[d]^p \\
				& \mathscr{X} \ar[r]^f & \mathscr{Y}
			}
		\end{align}
		where $X$ and $Y$ are algebraic spaces, the square is cartesian, g and p are smooth and surjective. If we further require that $X$ and $Y$ are schemes, then we call this \textit{a chart of f by schemes}\\
		Let $P$ be a property of morphisms of schemes which is stable and local on domain with respect to the smooth topology. We say that a morphism $f:\mathscr{X}\rightarrow \mathscr{Y}$ has property P if there exists a chart for $f$ by schemes such that the morphism $h$ in diagram \ref{morphism_stack_diagram} has property $P$.
	\end{definition}  
	
	\begin{Prpn}
		Let $f:\mathscr{X}\rightarrow \mathscr{Y}$ be a morphism of algebraic stacks over $S$ and the diagram below be \textit{a chart for f by schemes}
		\begin{center}
			$\xymatrix{
				X \ar[r]^g \ar[rd]^q \ar@/^1.5pc/[rr]^h & \mathscr{X}' \ar[r]^{f'} \ar@{->}[d]^{p'} & Y \ar@{->}[d]^p \\
				& \mathscr{X} \ar[r]^f & \mathscr{Y}
			}$
		\end{center}
		then the diagram obtained after applying the functor $J^{n}$, represented below is also \textit{a chart for f by schemes} 
		
		\begin{align}\label{Propofmorph}
			\xymatrix{
				J^{n}X \ar[r]^{J^{n}(g)} \ar[rd]_{J^{n}(q)}\ar@/^1.8pc/[rr]^{J^{n}(h)} & J^{n}\mathscr{X}' \ar[r]^{J^{n}(f')} \ar@{->}[d]^{J^{n}(p')} & J^{n}Y \ar@{->}[d]^{J^{n}(p)} \\
				& J^{n}\mathscr{X} \ar[r]^{J^{n}(f)} & J^{n}\mathscr{Y}
			}
		\end{align}

	\end{Prpn}  
	
	\begin{proof}
		Note that $J^{n}X$ and $J^{n}Y$ are schemes, using Proposition \ref{lem:surjection},  $J^{n}X\xrightarrow{{J^{n}(g)}} J^{n}\mathscr{X}'$ and $J^{n}Y\xrightarrow{J^{n}(p)}J^{n}\mathscr{Y}$ are smooth and surjective, also the induced diagram is cartesian since $J^{n}$ preserves fiber products. Hence $J^{n}$ preserves all the properties required for \textit{a chart for f by schemes}. 
	\end{proof}	
	
	\begin{Cor}
		The following properties morphisms which are stable and local on domain in $(Sch/S)_{\acute{e}t}$ of algebraic stacks are preserved under the $J^{n}$ functor.
		\begin{enumerate}[label=(\roman*)]
			\item locally of finite type, locally of finite presentation.
			\item étale, smooth
		\end{enumerate}
	\end{Cor}
	
	\begin{proof}
		Note that from diagram \ref{Propofmorph}, $J^{n}(f)$ satisfy a property $P$ which is stable and local on the domain if $J^{n}(h)$ has the property $P$. Now the above listed properties which are stable and local on the domain are preserved by $J^{n}$ for schemes. That is \cite[Proposition ~11.1]{article3} for étale and smooth morphisms, \cite[Proposition ~11.3 (c),(d)]{article3} for locally of finite type and locally of finite presentation.
	\end{proof}
	\begin{rmk}
		We have previously proved the jet space functor preserves properties like étale and smoothness when the given the functor is representable. In the case above, there is no representability condition imposed on the morphism between stacks.
	\end{rmk}
	
\printbibliography
	\nocite{*}
	
\end{document}